\documentclass[a4paper,11pt]{article}
\usepackage{amsmath,amssymb,amsthm}
\usepackage{cases}
\usepackage{a4wide}
\usepackage{color,graphicx}
\usepackage[hang,small,bf]{caption}
\usepackage[subrefformat=parens]{subcaption}
\usepackage{here}
\usepackage{comment}
\usepackage{ulem}

\numberwithin{equation}{section}
\newtheorem{thm}{Theorem}[section]
\newtheorem{prop}[thm]{Proposition}
\newtheorem{lem}[thm]{Lemma}
\newtheorem{rem}[thm]{Remark}

\newtheorem{asm}{Assumption}

\newcommand{\dL}[1]{\,\mathrm{d}\mathcal{L}^{#1}}
\newcommand{\dH}[1]{\,\mathrm{d}\mathcal{H}^{#1}}
\newcommand{\dt}{\,\mathrm{d}t}
\newcommand{\dd}[1]{\frac{\rm d}{{\rm d}#1}}
\newcommand{\ddt}{\dd{t}}

\newcommand{\bR}{{\mathbb R}}

\newcommand{\pOmega}[0]{\partial\Omega}
\newcommand{\closure}[1]{\overline{#1}}

\newcommand{\curve}[2]{\Gamma^{#1}_{#2}}

\newcommand{\chemicalVec}[2]{\bv{w}^{#1}_{#2}}
\newcommand{\chemicalDisc}[2]{{W^{#1}_{#2}}}
\newcommand{\chemicalDiscVec}[2]{{\bv{W}^{#1}_{#2}}}
\newcommand{\identity}[0]{\vec{\operatorname{Id}}}
\newcommand{\X}[2]{\vec{X}^{#1}_{#2}}
\newcommand{\innerproduct}[4]{\left<#1,#2\right>^{#3}_{#4}}

\newcommand{\normal}[3]{\vec{\nu}^{#1}_{#2,#3}}
\newcommand{\normalnoindex}[2]{\vec{\nu}^{#1}_{#2}}
\newcommand{\normalcontainer}[0]{\vec{\nu}_\Omega}

\newcommand{\conormal}[1]{\vec{\mu}_{#1}}
\newcommand{\vertexnormalnoindex}[2]{\vec{\nu}^{#1}_{#2}}

\newcommand{\weightnormal}[3]{\vec{\omega}^{#1}_{#2,#3}}
\newcommand{\weightnormalnoindex}[2]{\vec{\omega}^{#1}_{#2}}
\newcommand{\spacemeshpointsnum}[0]{K^m_{\Omega}}
\newcommand{\sgrad}[0]{\nabla_s} 
\newcommand{\slaplace}[0]{\Delta_s} 
\newcommand{\vertex}[3]{\vec{q}^{#1}_{#2,#3}}
\newcommand{\edge}[3]{e^{#1}_{#2,#3}}

\newcommand{\Glinear}[2]{\pi^{h}_{#2}\left[\frac{\X{#1+1}{#2} - \X{#1}{#2}}{\tau_{#1}}\cdot\weightnormalnoindex{#1}{#2}\right]}

\newcommand{\basisbulk}[2]{\Psi^{#1}_{#2}}
\newcommand{\metric}[1]{\uuline{#1}}
\newcommand{\zerovec}[0]{\bvzero}
\newcommand{\zerosumfunc}[0]{T\widetilde{\Sigma}}
\newcommand{\zerosumspace}[0]{T\Sigma}

\newcommand{\femspacebulk}[0]{{S}^m}
\newcommand{\femspacebulkvector}[0]{\bv {S}^m}
\newcommand{\femspacebulkvectorsum}[0]{\bv {S}^m_\Sigma}
\newcommand{\femspacecurvevector}[1]{{\underline{V}^h_{#1}}}
\newcommand{\femspacegamma}[0]{V(\Gamma^m)}
\newcommand{\femspacegammavector}[0]{\underline{V}(\Gamma^m)}
\newcommand{\femspacecurve}[1]{{V^h_{#1}}}
\newcommand{\basiscurve}[2]{{\Phi_{{#1},{#2}}}}

\newcommand{\triplejunction}[1]{\mathcal{T}_{#1}}
\newcommand{\curveindex}[2]{c^{#1}_{#2}}
\newcommand{\curveindexrho}[2]{\rho^{#1}_{#2}}

\newcommand{\jump}[1]{\left[{#1}\right]}

\newcommand{\naturalset}[1]{\mathbb{N}_{\leq {#1}}}
\newcommand{\phasecharacter}[0]{\bv \chi}

\newcommand{\dcmap}[0]{\mathcal O}

\newcommand{\constalpha}[1]{\frac{\frac{1}{R_1{#1}} + \frac{1}{R_2{#1}} + \frac{1}{R_3{#1}}}{2\log{\frac{R_3{#1}}{R_1{#1}}}}}

\def\bv#1{\mbox{\boldmath{$#1$}}}    

\def\bvzero{{\bf 0}}
\def\bvone{{\bf 1}}
\newcommand{\bigchi}{\ensuremath{\mathrm{\mathcal{X}}}}
\newcommand{\charfcn}[1]{\bigchi_{#1}}

\title{A structure-preserving finite element method for the multi-phase Mullins--Sekerka problem with triple junctions}
\author{Tokuhiro Eto\thanks{Graduate School of Mathematical Sciences, The University of Tokyo, Komaba 3-8-1, Meguro, Tokyo 153-8914, Japan. E-mail: tokuhiro\_eto@yahoo.co.jp}\and Harald Garcke\thanks{Fakult\"{a}t f\"{u}r Mathematik, Universit\"{a}t Regensburg, 93040 Regensburg, Germany. E-mail: harald.garcke@ur.de}\and Robert N\"{u}rnberg\thanks{Dipartimento di Mathematica, Universit\`{a} di Trento, 38123 Trento, Italy. E-mail: robert.nurnberg@unitn.it}}
\date{}

\begin{document}

\maketitle

\begin{abstract}
We consider a sharp interface formulation for the multi-phase Mullins--Sekerka flow. The flow is characterized by a network of curves evolving such that the total
surface energy of the curves is reduced, while the areas of the enclosed phases are
conserved. Making use of a variational formulation, we introduce a fully
discrete finite element method. Our discretization features a parametric
approximation of the moving interfaces that is independent of the
discretization used for the equations in the bulk.
The scheme can be shown to be unconditionally stable and to satisfy
an exact volume conservation property. Moreover, an inherent tangential 
velocity for the vertices on the discrete curves leads to asymptotically
equidistributed vertices, meaning no remeshing is necessary in practice.
Several numerical examples, including a convergence experiment for the three-phase
Mullins--Sekerka flow, demonstrate the capabilities of the introduced method.
\end{abstract}


\section{Introduction}
In this paper, we consider the problem of networks of curves moving under the
multi-phase Mullins--Sekerka flow, see, e.g., \cite{BronsardGarckeStoth1998}. 
These networks feature triple junctions, at which certain balance laws need to
hold.
The network depicted in Figure~\ref{fig:problem_setting} consists of three
time dependent curves $\curve{}{i}(t)$, $i=1,2,3$, that meet
at two triple junction points $\mathcal{T}_1(t)$ and $\mathcal{T}_2(t)$.
We assume that the curve network $\curve{}{}(t) := \cup_{i=1}^3\curve{}{i}(t)$ 
lies in a domain $\Omega\subset\mathbb{R}^2$
for all times $t \in [0, T ]$, and it partitions
$\Omega$ into the three subdomains $\Omega_j(t)$, $j=1,2,3$. 
The three domains correspond to different phases in the multi-component system.
The evolution of the interfaces $\curve{}{1}(t),\curve{}{2}(t),$ and $\curve{}{3}(t)$ is driven by diffusion.
As in \cite{BronsardGarckeStoth1998}, given a time $T > 0$ and the hyperplane 
$T\Sigma := \{\bv u\in\mathbb{R}^3\mid \sum_{j=1}^3u_j = 0\}$,
we introduce a vector of chemical potentials $\chemicalVec{}{}:(0,T)\times\Omega\to T\Sigma$ which fulfills 
the quasi-static diffusion equation for $j=1,2,3$ and $t\in(0,T)$,
\begin{subequations} \label{eq:SF}
\begin{equation}\label{eq:StrongForm_LaplaceEq}
    \Delta\chemicalVec{}{} = \zerovec\ \ \mbox{in}\ \ \Omega_j(t)
\end{equation}
together with
\begin{equation}\label{eq:StrongForm_Neumann}
    \partial_{\normalcontainer}\chemicalVec{}{} = \zerovec\ \ \mbox{on}\ \ \pOmega,
\end{equation}
where $\normalcontainer$ denotes the {outer} unit normal vector to $\pOmega$. 
\begin{figure}[H]
    \centering
    \includegraphics[keepaspectratio, scale=0.25]{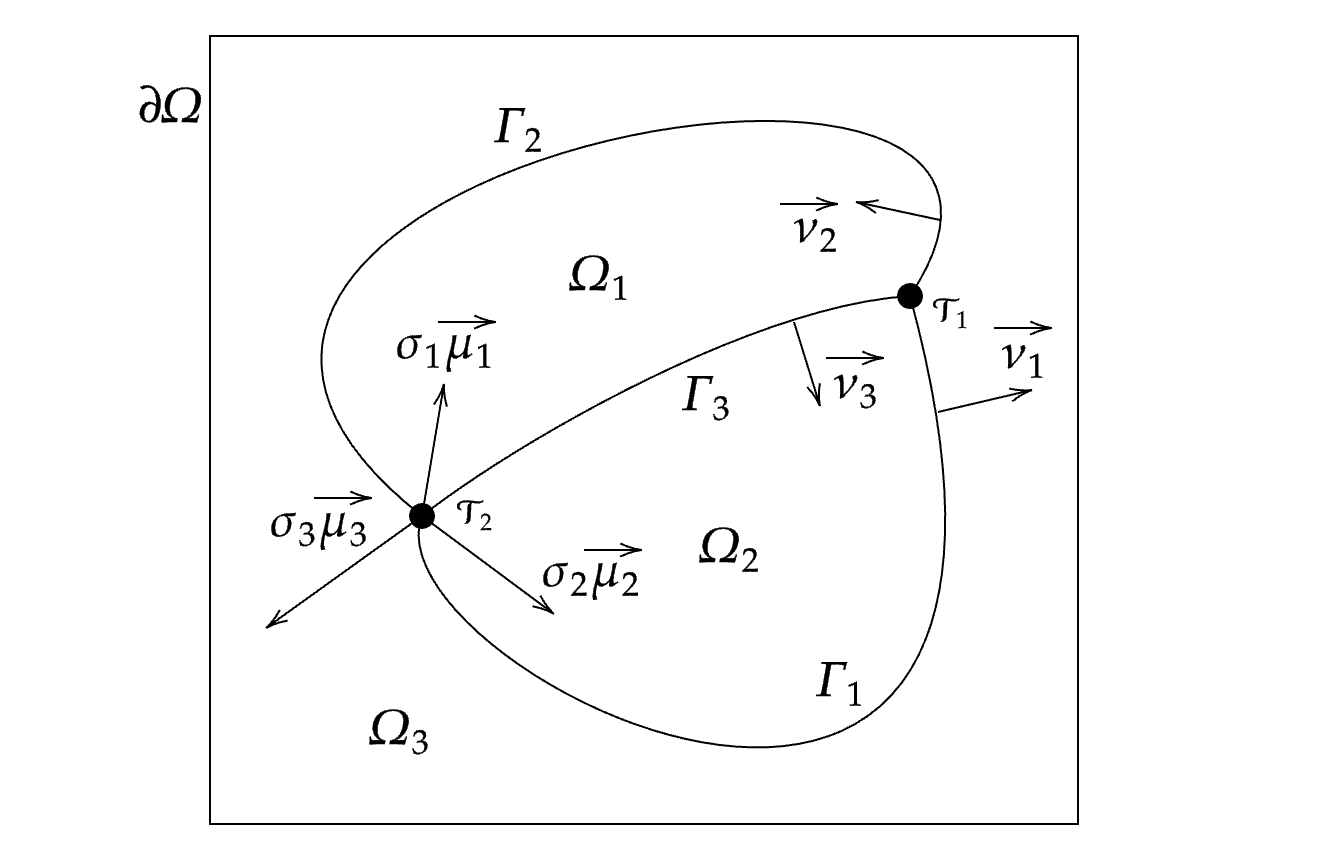}
    \caption{Three open curves with two triple junctions.}\label{fig:problem_setting}
\end{figure}
To close the system, we need boundary conditions on $\curve{}{}(t)$
and on $\mathcal{T}_1(t)$ and $\mathcal{T}_2(t)$.
These boundary conditions are given by a 
Stefan-type kinetic condition and the Gibbs--Thomson law on the moving
interfaces, and Young's law at the triple junctions, see 
\cite{BronsardGarckeStoth1998}.
The kinetic condition reads
\begin{equation}\label{eq:StrongForm_MotionLaw}
    [\nabla\chemicalVec{}{}]\vec{\nu} = -V[\phasecharacter]\ \ \mbox{on}\ \ \curve{}{}(t),
\end{equation}
where $\phasecharacter = (\chi_1,\chi_2,\chi_3)^T$ denotes the vector which consists of the characteristic functions $\chi_j = \charfcn{\Omega_j(t)}$ of 
$\Omega_j(t)$, $\vec{\nu}$ is the unit normal vector on $\curve{}{}(t)$,
and $V$ is the velocity of $\Gamma(t)$ in the direction of $\vec{\nu}$.
We write $\vec\nu = \sum_{i=1}^3 \vec\nu_i\charfcn{\Gamma_i(t)}$
and use this convention for quantities defined on $\Gamma(t)$ throughout the
paper.
The orientation of the three normal vectors is shown in 
Figure~\ref{fig:problem_setting}.
In addition, the quantity $[q]$ represents the jump of $q$ across $\Gamma(t)$ in the direction of $\vec{\nu}$ defined by $[q](\vec x) := \lim_{\varepsilon\searrow 0}\{q(\vec x + \varepsilon\vec{\nu}) - q(\vec x - \varepsilon\vec{\nu})\}$. 
Furthermore, the Gibbs--Thomson equations can be written as
\begin{equation}\label{eq:StrongForm_GibbsThomsonLaw}
    \chemicalVec{}{}\cdot[\phasecharacter] = \sigma\varkappa\ \ \mbox{on}\ \ \curve{}{}(t),
\end{equation}
where $\varkappa$ denotes the curvature of $\Gamma(t)$ (well-defined
on the interiors of $\Gamma_i(t)$), and
$\sigma= \sum_{i=1}^3 \sigma_i \charfcn{\Gamma_i(t)}$ is a surface 
tension coefficient on $\Gamma(t)$. Our sign convention is such that unit
circles have curvature $\varkappa=-1$, which is different to the one used in
\cite{BronsardGarckeStoth1998}.
Finally, denoting by $\conormal{i}$, the outer unit co-normal to $\curve{}{i}(t)$, we further require Young's law,
which is a balance of force condition at the triple junction as follows:
\begin{equation}\label{eq:YoungLaw}
    \sum_{i=1}^3 \sigma_i\conormal{i} = \vec 0\ \ \mbox{on}\ \ \partial\curve{}{1}(t)\cap\partial\curve{}{2}(t)\cap\partial\curve{}{3}(t).
\end{equation}
\end{subequations}
In order to be able to fulfill this condition, we require
$\sigma_1 \leq \sigma_2 + \sigma_3$, $\sigma_2 \leq \sigma_1 + \sigma_3$ and
$\sigma_3 \leq \sigma_1 + \sigma_2$.
It can be shown that solutions to
\eqref{eq:SF} reduce the weighted length $\sum_{i=1}^3 \sigma_i |\Gamma_i(t)$
of the curve network, while conserving the areas of the subdomains
$\Omega_1(t)$ and $\Omega_2(t)$ (and hence trivially also of $\Omega_3(t)$),
see Section~\ref{sec:MathematicalProperties} for the precise details.
Our aim in this paper is to introduce a numerical method that preserves these
two properties on the discrete level.

Prescribing an initial condition $\curve{}{}(0) = \curve{0}{}$ for the interface, we altogether obtain the following system:
\begin{equation}\label{eq:StrongForm}
    \begin{cases}
        \Delta\bv w = \zerovec\ \ \mbox{in}\ \ (0,T]\times(\Omega\backslash\curve{}{}(t)),\\
        \bv w\cdot[\phasecharacter] = \sigma \varkappa\ \ \mbox{on}\ \ \curve{}{}(t),\ t\in[0,T],\\
        [\nabla\bv w]\,\vec{\nu} = -V[\phasecharacter]\ \ \mbox{on}\ \ \curve{}{}(t),\ t\in[0,T],\\
        \partial_{\normalcontainer}\bv w = \zerovec\ \ \mbox{on}\ \ \pOmega,\\
        \sum_{i=1}^3\sigma_i\vec{\mu}_i = \vec 0\ \ \mbox{on}\ \ \partial\curve{}{1}(t)\cap\partial\curve{}{2}(t)\cap\partial\curve{}{3}(t),\ t\in[0,T],\\
        \curve{}{}(0) = \Gamma^0.
    \end{cases}
\end{equation}
\newline

The system \eqref{eq:StrongForm} at present is written for the setup from
Figure~\ref{fig:problem_setting}, i.e.\ a network of three curves, meeting at
two triple junctions and partitioning $\Omega$ into three phases.
We will later generalize this to an arbitrary network of curves. The simplest
case is given by a single closed curve that partitions the domain into two
phases. Then we obtain the classical
two phase Mullins--Sekerka problem, 
see \cite{ChenHongYi1996} and the references given below. Indeed, let 
$(\bv w,\{\curve{}{}(t)\}_{0\leq t\leq T})$, with $\bv w = (w_1,w_2)^T$, 
be a solution to the corresponding problem \eqref{eq:StrongForm}
with $\sigma=1$, and let
$\normalnoindex{}{}$ point into $\Omega_2(t)$, the interior domain of $\curve{}{}(t)$.
Then we have 
$[\phasecharacter] = {(-1,1)}^T$ on $\curve{}{}(t)$, and it holds that 
$w = {w_2 - w_1}$ is a solution to the system
\begin{equation}\label{eq:2-StrongForm}
    \begin{cases}
        \Delta w = 0\ \ \mbox{in}\ \ [0,T]\times(\Omega\backslash\curve{}{}(t)),\\
        w = \varkappa\ \ \mbox{on}\ \ \curve{}{}(t),\ \ t\in[0,T], \\
        \frac{1}{2}[\nabla w] \cdot \vec\nu = -V \ \ \mbox{on}\ \ \curve{}{}(t),\ \ t\in[0,T],\\
        \partial_{\normalcontainer}w = 0\ \ \mbox{on}\ \ \pOmega.
    \end{cases}
\end{equation}

The multi-phase Mullins--Sekerka problem \eqref{eq:StrongForm} arises naturally
as the sharp interface limit of a nondegenerate multi-component Cahn--Hilliard
equation.
Let $\Sigma:=\{\bv u\in\mathbb{R}^3\mid\sum_{j=1}^3u_j = 1\}$ and 
let $\zerosumfunc$ be the family of all functions $\chemicalVec{}{}:\Omega\to\mathbb{R}^3$ 
such that $\operatorname{Im}{(\chemicalVec{}{})}\subset\zerosumspace$. 
Let $\psi:\mathbb{R}^3\to\mathbb{R}$ be a potential whose restriction to $\Sigma$
has exactly three distinct and strict global
minima, say $\bv p_i\in\Sigma$, $i=1,2,3$, with $\psi(\bv p_1) = \psi(\bv p_2) = \psi(\bv p_3)$.
Let $F:\mathbb{R}^3\to\mathbb{R}^3$ be the projection of $\nabla {\psi}$ onto $\zerosumfunc$.
According to \cite[\S2]{BronsardGarckeStoth1998}, 
the system $\eqref{eq:StrongForm}$ is derived as the limit with $\varepsilon\to 0$ of
a chemical system consisting of three species governed by 
the vector-valued Cahn--Hilliard equation whose form reads as:
\begin{equation}\label{eq:CahnHilliardVV}
    \begin{cases}
        \partial_t \bv {u} = \Delta \bv{w}\ \ \mbox{in}\ \ \Omega\times(0,T],\\
        \bv{w} = -\varepsilon\Delta \bv{u} + \frac{1}{\varepsilon}F(\bv{u})\ \ \mbox{in}\ \ \Omega\times(0,T],\\
        \partial_{\normalcontainer} \bv{u} = \partial_{\normalcontainer} \bv{w} = 0\ \ \pOmega\times[0,T],\\
        \bv{u}(0,\cdot) = \bv{u}_0\ \ \mbox{in}\ \ \Omega,
    \end{cases}
\end{equation}
where $\bv u_0:\Omega\to\Sigma$ denotes the initial distribution of each 
component and
$\bv{u}:\Omega\times[0,T]\to\Sigma$ and $\chemicalVec{}{}:\Omega\times[0,T]\to T\Sigma$ 
indicate the concentration and the chemical potential of each component in
time, respectively.
A distributional solution concept to $\eqref{eq:StrongForm}$ was proposed, and its existence 
was established via an implicit time discretization and
under the assumption that no interfacial energy is lost in the limit in the time discretization (see \cite[Definition~4.1, Theorem~5.8]{BronsardGarckeStoth1998}).
See \cite{GarckeSturzenhecker1998} for a related work which treated the case without triple junctions and with a 
driving force.\newline

Compared to the multi-phase Mullins--Sekerka problem, 
the binary case, namely {the} two-phase {case}, has been well studied so far. For classical solutions, 
Chen et al. \cite{ChenHongYi1996} showed the existence 
of a classical solution to the Mullins--Sekerka problem local-in-time in {the} two-dimensional case, whereas
Escher and Simonett \cite{EScherSimonett} gave a similar result in the general dimensional case.
When it comes to the notion of weak solution{s}, Luckhaus and Sturzenhecker \cite{LuckhausSturzenhecker}
established the existence of weak solutions to $\eqref{eq:StrongForm}$ in a distributional sense.
Therein, the weak solution was obtained as a limit of a sequence of time discrete approximate solutions
under the no mass loss assumption. {The} time implicit scheme is {the basis of the approach} in \cite{BronsardGarckeStoth1998}.
After that, R\"{o}ger \cite{Roger2005} removed the technical assumption {of no mass loss} in the case when the Dirichlet--Neumann
boundary condition is imposed by using geometric measure theory. 
Recently, researches which treat the boundary contact case
gradually appear. Garcke and Rauchecker \cite{GarckeRauchecker2022} {considered a} stability analysis
in a curved domain in $\mathbb{R}^2$ via a {linearization approach}.
Hensel and Stinson \cite{hensel2022weak} 
proposed a varifold solution to $\eqref{eq:StrongForm}$ by starting from the energy dissipation property. 
For a gradient flow aspect of the Mullins--Sekerka flow, see e.g., \cite[\S3.2]{Serfaty2010}.
\newline

The numerical scheme that we propose in this paper
is based on the BGN method, a parametric finite element method
that allows the variational treatment of triple junctions and was first
introduced by Barrett, Garcke, and N\"{u}rnberg in 
\cite{BarrettGarckeRobert2007,BarrettGarckeNuernbergSurfaceDiffusion2007}. 
We also refer to the review article
\cite{BarrettGarckeRobertBook2020} for more details on the BGN method, including in the context 
of the standard Mullins--Sekerka problem \eqref{eq:2-StrongForm}. 
Alternative front-tracking methods for geometric flows of curve networks have
been considered in, e.g.,
\cite{BronsardW95,Thaddey99,Neubauer2002,PanW08,PozziS21}.

Let us briefly review numerical methods being available in the literature for
the Mullins--Sekerka problem {and for} its diffuse interface model, 
the multi-component Cahn--Hilliard equation \eqref{eq:CahnHilliardVV}.
To the best of our knowledge, there are presently no sharp interface
methods for the numerical approximation of 
the multi-phase Mullins--Sekerka problem.
For the boundary integral method for the two-phase case, we refer the reader to \cite{BatesBrown,Bates1995ANS,ChenKublikTsai2017,Mayer2000,ZhuChenHou}.
A level set formulation of moving boundaries together with the finite difference method was proposed in \cite{ChenMettimanOsherSmereka1997}.
For an implementation of the method of fundamental solutions for the Mullins--Sekerka problem in $2$D, see \cite{Eto2024}.
For the parametric finite element method in general dimensions, see 
\cite{BarrettGarckeRobert2010,Nurnberg202203}.
Numerical analysis of the scalar Cahn--Hilliard equation is dealt with in the works \cite{BarrettBlowey1995,BarrettBloweyGarcke1999,BloweyElliott1992,ElliottFrench1987}.
Feng and Prohl \cite{FengProhl2004} proposed a mixed fully discrete finite element method for the Cahn--Hilliard equation 
and provided error estimates between the solution of the Mullins--Sekerka problem and the approximate solution of the Cahn--Hilliard equation which are computed by their scheme.
The established error bounds yielded a convergence result in \cite{FengProhl2005}.
Aside from the sharp interface model, the Cahn--Hilliard equation for the multi-component case
has been computed in several works, see {\cite{BarrettBloweyGarcke2001, BloweyCopettiElliott1996,Eyre1993,LiChoiKim2016,Nurnberg09}}.
The multi-component Cahn--Hilliard equation on surfaces has recently been
considered in \cite{LiLiuXiaHeLi2022}.
\newline

This paper is organized as follows.
In the first part, we focus on the three-phase case, as outlined in the
introduction.
In Section~\ref{sec:MathematicalProperties}, we show a curve-shortening property and an area-preserving property of strong solutions to the system.
In Section~\ref{sec:WeakFormulation}, we introduce a weak formulation of the system, which in Section~\ref{sec:FEM_approx} will then be
approximated with the help of an unfitted parametric finite element method. 
The scheme, which is linear, can be shown to be unconditionally stable.
{Section~\ref{sec:MatrixForm} is devoted to discussing solution methods for
the linear systems that arise at each time level.}
In Section~\ref{sec:ImprovementOfTheScheme}{,} we adapt our approximation to allow
for an exact area-preservation property on the fully discrete level,
leading to a nonlinear scheme.
Section~\ref{sec:GeneralizationToMulti} is devoted to generalizations 
of the problem formulation and our numerical approximation to the
general multi-phase case.
Finally, we will show several results of numerical computation in 
Section~\ref{sec:NumericalResults}, including convergence experiments for 
a constructed solution in the three-phase case.

\section{Mathematical properties}\label{sec:MathematicalProperties}
In this section, we shall present two important properties 
of strong solutions to $\eqref{eq:StrongForm}$.

\begin{prop}[Curve shortening property of strong solutions]\label{prop:CurveShorteningStrong}
        Assume that $(\bv w,\{\curve{}{}(t)$\linebreak $\}_{0\leq t\leq T})$
        is a solution to $\eqref{eq:StrongForm}$.
        Then it holds that
\[
\ddt |\Gamma(t)|_\sigma + \int_{\Omega}|\nabla \bv w|^2\dL{2} = 0,
\]
where we have defined the weighted length $|\Gamma(t)|_\sigma = \sum_{i=1}^3 \sigma_i|\curve{}{i}(t)| = \sum_{i=1}^3 \sigma_i \int_{\curve{}{i}(t)} 1 \dH{1}$. Here $\dL{2}$ and $\dH{1}$ refer to integration with respect to the Lebesgue measure and the one-dimensional Hausdorff measure in $\bR^2$,
respectively.
\end{prop}
\begin{proof}
        Since the curvature is the first variation of the curves' length in the normal direction, we deduce, on recalling \eqref{eq:StrongForm_GibbsThomsonLaw} 
and \eqref{eq:StrongForm_MotionLaw}, that
        \begin{align}\label{eq:CSProp_1}
\ddt\sum_{i=1}^3\sigma_i|\curve{}{i}(t)| & = -\sum_{i=1}^3 \sigma_i\int_{\curve{}{i}(t)} \varkappa_i V_i\dH{1}
= - \int_{\curve{}{}(t)} \sigma \varkappa V\dH{1} 
= - \int_{\curve{}{}(t)} \chemicalVec{}{}\cdot[\phasecharacter] V\dH{1} 
\nonumber \\ & 
= \int_{\curve{}{}(t)} \chemicalVec{}{}\cdot
[\nabla\chemicalVec{}{}]\vec{\nu} \dH{1}
= \sum_{j=1}^3 \int_{\curve{}{}(t)} w_j\, [\nabla w_j]\cdot\vec{\nu} \dH{1}.
        \end{align}
Now for each $1 \leq j \leq 3$ 
it follows from integration by parts, \eqref{eq:StrongForm_LaplaceEq} and
\eqref{eq:StrongForm_Neumann} that
        \begin{equation}\label{eq:CSProp_4}
\int_{\Omega} |\nabla w_j|^2\dL{2} = 
\sum_{k=1}^3 \int_{\Omega_k(t)} |\nabla w_j|^2\dL{2} = 
- \int_{\curve{}{}(t)} w_j\, [\nabla w_j]\cdot\vec{\nu} \dH{1}.
        \end{equation}
Summing \eqref{eq:CSProp_4} for $j=1,2,3$ 
and combining with \eqref{eq:CSProp_1} yields the desired result.
\end{proof}

\begin{prop}[Area preserving property of strong solutions]\label{prop:AreaPreservingStrong}
Let $\Omega_1(t),\Omega_2(t),$ and $\Omega_3(t)$ be the domains bounded by $\curve{}{2}(t)\cup\curve{}{3}(t)$, $\curve{}{1}(t)\cup\curve{}{3}(t)$, and $\curve{}{1}(t)\cup\curve{}{2}(t)\cup\pOmega$, respectively, see Figure~\ref{fig:problem_setting}. Then for any solution to $\eqref{eq:StrongForm}$, it holds that
\[
\ddt |\Omega_j(t)| = 0 \,,\qquad j = 1,2,3.
\]
\end{prop}
\begin{proof}
We first deduce from the motion law $\eqref{eq:StrongForm_MotionLaw}$ that
\begin{flalign}\label{eq:APProp_1}
    \begin{split}
        \ddt|\Omega_1(t)| 
        &= \int_{\curve{}{2}(t)}V_2\dH{1} + \int_{\curve{}{3}(t)}-V_3\dH{1} \\
        &= \int_{\curve{}{2}(t)}-\jump{\chi_1}\jump{\nabla w_1}\cdot\normalnoindex{}{}\dH{1} + \int_{\curve{}{3}(t)}\jump{\chi_1}\jump{\nabla w_1}\cdot\normalnoindex{}{}\dH{1} \\
        &=\int_{\curve{}{2}(t)}-\jump{\nabla w_1}\cdot\normalnoindex{}{}\dH{1} + \int_{\curve{}{3}(t)}-\jump{\nabla w_1}\cdot\normalnoindex{}{}\dH{1} \\
        &= \int_{\pOmega_1(t)\cap\curve{}{}(t)}-\jump{\nabla w_1}\cdot\normalnoindex{}{}\dH{1}= \int_{\curve{}{}(t)}-\jump{\nabla w_1}\cdot\normalnoindex{}{}\dH{1}.
    \end{split}
\end{flalign}
Here, we note that $\nabla w_1$ does not jump over $\curve{}{1}(t)$ to derive the last equality.
Meanwhile, integration by parts with $\eqref{eq:StrongForm_LaplaceEq}$ and $\eqref{eq:StrongForm_Neumann}$ shows
\begin{equation}\label{eq:APProp_2}
    0 = \int_{\Omega(t)}\Delta w_1\dL{2} = \int_{\Omega_1(t)}\Delta w_1\dL{2} + \int_{\Omega\backslash\closure{\Omega_1(t)}}\Delta w_1\dL{2} = \int_{\curve{}{}(t)}-\jump{\nabla w_1}\cdot\normalnoindex{}{}\dH{1}.
\end{equation}
Combining $\eqref{eq:APProp_1}$ and $\eqref{eq:APProp_2}$ gives the assertion.
The cases for $j=2,3$ can be shown in the same manner.
\end{proof}

\section{Weak formulation}\label{sec:WeakFormulation}
Let us derive a weak formulation {for} $\eqref{eq:StrongForm}$.
In the sequel, we often abbreviate $\Omega_j(t)$ and $\curve{}{i}(t)$ as $\Omega_j$ and $\curve{}{i}$, for simplicity.
Now $1\leq j\leq 3$.
Then, testing the first equation of $\eqref{eq:StrongForm}$ with $\varphi\in H^1(\Omega)$, 
we deduce similarly to $\eqref{eq:APProp_2}$ that
\begin{equation*}
    0 = \int_{\Omega}\Delta w_j\varphi\dL{2} = \sum_{k=1}^3\int_{\Omega_k}\Delta w_j\varphi\dL{2} = \sum_{k=1}^3\int_{\curve{}{}(t)}-\jump{\nabla w_j}\cdot\normalnoindex{}{}\varphi\dH{1} - \int_\Omega \nabla w_j\cdot\nabla\varphi\dL{2}.
\end{equation*}
Then, we obtain from the second equality of $\eqref{eq:APProp_1}$ that
\begin{equation}\label{eq:WeakForm1}
    \begin{cases}
        \int_{\curve{}{2}}V_2\varphi\dH{1} - \int_{\curve{}{3}}V_3\varphi\dH{1} - \int_\Omega \nabla w_1\cdot\nabla\varphi\dL{2} = 0,\\
        \int_{\curve{}{3}}V_3\varphi\dH{1} - \int_{\curve{}{1}}V_1\varphi\dH{1} - 
        \int_{\Omega}\nabla w_2\cdot\nabla\varphi\dL{2} = 0, \\
        \int_{\curve{}{1}}V_1\varphi\dH{1} - \int_{\curve{}{2}}V_2\varphi\dH{1} - 
        \int_{\Omega}\nabla w_3\cdot\nabla\varphi\dL{2} = 0.
    \end{cases}
\end{equation}
Hence{,} $\bv w \in \zerosumfunc \cap [H^1(\Omega)]^3$ is such that
\begin{equation} \label{eq:WeakForm1a}
\int_\Omega \nabla \bv w : \nabla \bv \varphi -
\sum_{(i,j,k)\in\Lambda} \int_{\curve{}{i}}V_i(\varphi_j - \varphi_k)\dH{1} = 0 \qquad
\forall\ \bv \varphi \in \zerosumfunc \cap [H^1(\Omega)]^3{,}
\end{equation}
where
\begin{equation}\label{eq:IndexTriplet}
    \Lambda := \{(1,3,2),(2,1,3),(3,2,1)\}.
\end{equation}

The Gibbs--Thomson relation $\eqref{eq:StrongForm_GibbsThomsonLaw}$ is encoded as follows:
\begin{equation}\label{eq:WeakForm2}
    \begin{cases}
        \int_{\curve{}{1}}(w_2 - w_3 + \sigma_1 \varkappa_1)\xi\dH{1} = 0\ \ \forall \xi\in L^2(\curve{}{1}), \\
        \int_{\curve{}{2}}(w_3 - w_1 + \sigma_2 \varkappa_2)\xi\dH{1} = 0\ \ \forall\xi\in L^2(\curve{}{2}), \\
        \int_{\curve{}{3}}(w_1 - w_2 + \sigma_3 \varkappa_3)\xi\dH{1} = 0\ \ \forall\xi\in L^2(\curve{}{3}).
    \end{cases}
\end{equation}
Finally, we give a weak formulation of the weighted curvature vector
$\varkappa_\sigma := \sigma\varkappa$, which means that $\varkappa_\sigma = \sigma_i\varkappa_i$ on $\curve{}{i}(t)$ for $i=1,2,3$.
Let $\vec{\operatorname{Id}}$ denote the identity map in $\mathbb{R}^2$. 
Then, it holds that 
$\sigma\slaplace\vec{\operatorname{Id}} = \varkappa_\sigma\vec{\nu}$ on $\Gamma$, see \cite{Dziuk88}.
Take a test function $\vec{\eta}\in H^1(\Gamma;\mathbb{R}^2)$ with 
$\vec{\eta}\!\mid_{\curve{}{1}} = \vec{\eta}\!\mid_{\curve{}{2}} = \vec{\eta}\!\mid_{\curve{}{3}}$ 
on $\partial\curve{}{1}\cap\partial\curve{}{2}\cap\partial\curve{}{3}$. 
Then applying integration by parts gives
\begin{align}\label{eq:WeakForm3}
    \int_{\Gamma}\varkappa_\sigma\vec{\nu}\cdot\vec{\eta}\dH{1} & = \int_\Gamma\sigma\slaplace\identity\cdot\vec{\eta}\dH{1} = \sum_{i=1}^3\int_{\curve{}{i}}\sigma_i\slaplace\identity\cdot\vec{\eta}\dH{1} \nonumber \\
& = \sum_{i=1}^3\left(\sigma_i \int_{\partial\curve{}{i}}(\sgrad\identity\ \vec{\eta})\cdot\vec{\mu}_i\dH{0} - \sigma_i\int_{\curve{}{i}}\sgrad\identity\cdot\sgrad\vec{\eta}\dH{1}\right) \nonumber \\
& = \sum_{j=1}^2\left\{\left(\sum_{i=1}^3\sigma_i\vec{\mu}_i\right)\cdot(\sgrad\identity\ \vec{\eta})\right\}\!\mid_{\mathcal{T}_j} - \sum_{i=1}^3\sigma_i\int_{\curve{}{i}}\sgrad\identity\cdot\sgrad\vec{\eta}\dH{1} \nonumber \\
& = -\int_\Gamma\sigma\sgrad\identity\cdot\sgrad\vec{\eta}\dH{1}.
\end{align}
Here, we have used Young's law $\eqref{eq:YoungLaw}$ to get the last equality.
For later use, we define inner products on $\curve{}{}$ and $\Omega$ as follows:
\begin{equation*}
    \innerproduct{u}{v}{}{\Omega} := \int_\Omega u\,v\dL{2}\ \ \mbox{for}\ \ u,v\in L^2(\Omega),\qquad \innerproduct{u}{v}{}{\curve{}{}}:= \int_{\curve{}{}}u\,v\dH{1}\ \ \mbox{for}\ \ u,v\in L^2(\curve{}{}).
\end{equation*}
Let us summarize the weak formulation of the system $\eqref{eq:StrongForm}$ as follows, where we recall \eqref{eq:IndexTriplet}.
\newline\newline\noindent
\textbf{[Motion law]} For all $\bv\varphi\in\zerosumfunc\cap [H^1(\Omega)]^3$,
\begin{subequations} \label{eq:WFM}
\begin{equation}\label{eq:WFM-ML}
    \innerproduct{\nabla\bv w}{\nabla\bv \varphi}{}{\Omega} - \sum_{(i,j,k)\in\Lambda}\innerproduct{V_i}{\varphi_j-\varphi_k}{}{\curve{}{i}} = 0.
\end{equation}
\textbf{[Gibbs--Thomson law]} For all $\xi\in L^2(\curve{}{})$,
\begin{equation}\label{eq:WFM-GTL}
    \innerproduct{\varkappa_\sigma}{\xi}{}{\curve{}{}} - \sum_{(i,j,k)\in\Lambda}\innerproduct{w_j-w_k}{\xi_i}{}{\curve{}{i}} = 0.
\end{equation}
\textbf{[Curvature vector]} For all $\vec{\eta}\in H^1(\curve{}{};\mathbb{R}^2)$ with 
$\vec{\eta}_1 = \vec{\eta}_2 = \vec{\eta}_3$ 
on $\partial\curve{}{1}\cap\partial\curve{}{2}\cap\partial\curve{}{3}$,
\begin{equation}\label{eq:WFM-CV}
    \innerproduct{\varkappa_\sigma\normalnoindex{}{}}{\Vec{\eta}}{}{\curve{}{}} + \innerproduct{\sigma\sgrad\identity}{\sgrad\Vec{\eta}}{}{\curve{}{}} = 0.
\end{equation}
\end{subequations}

\section{Finite element approximation}\label{sec:FEM_approx}
To approximate the weak solution $(\bv{w},V,\varkappa_\sigma)$ of $\eqref{eq:StrongForm}$, we use ideas from \cite[(2.30a,b)]{BarrettGarckeRobert2007} and \cite[\S3]{Nurnberg202203}. Let the time interval $[0,T]$ be split into $M$ sub-intervals $[t_{m-1}, t_m]$ for each $m = 1,\cdots,M$ whose length are equal to $\tau_m$. 
Then, given a triplet of polygonal curves $\curve{0}{} = (\curve{0}{1},\curve{0}{2},\curve{0}{3})$, our aim is to find time discrete triplets $\curve{1}{},\cdots,\curve{M}{}$ governed by discrete analogues of $\eqref{eq:WeakForm1},\eqref{eq:WeakForm2}$ and $\eqref{eq:WeakForm3}$. 
For each $m\geq 0$ and $1\leq i\leq 3$, $\curve{m}{i} = \X{m}{i}(I)$ is parameterized by $I = [0,1]\ni\rho \mapsto \X{m}{i}(\rho)\in\mathbb{R}^2$ and $I$ is split into sub-intervals as $I = \cup_{j=1}^{N_i}[q_{i,j - 1}, q_{i,j}]$, where $N_i\in\mathbb{N}_{\geq2}$. 
Then, we note that
\begin{equation*}
    \innerproduct{u}{v}{}{\curve{m}{i}} = \int_{I}u\,v\,|(\X{m}{i})_s|\dL{1}.
\end{equation*}
Set $\vertex{m}{i}{j} := \X{m}{i}(q_{i,j})$. Let $\mathfrak{T}^m$ be a sequence of triangulations of $\closure{\Omega}$ and let
$\femspacebulk$ be the associated scalar- and vector-valued finite element 
spaces, namely
\begin{equation*}
    \femspacebulk := \left\{v\in C(\closure{\Omega})\mid v\!\mid_o\ \ \mbox{is affine}\ \ \forall o\in\mathfrak{T}^m\right\}, \quad 
\femspacebulkvector := [\femspacebulk]^3\quad\text{and}\quad
    \femspacebulkvectorsum := \femspacebulkvector \cap \zerosumfunc.
\end{equation*}
Let $\femspacecurve{i}$ be the set of all piecewise continuous functions on {$I$} which are affine on each {sub-interval $[q_{i,j-1},q_{i,j}]$},
and let $\pi^{{h}}_i:C({I})\to \femspacecurve{i}$ be the associated
standard interpolation operators {for $1\leq i\leq 3$}.
Similarly, $\femspacecurvevector{i}$ denotes the set of all vector valued functions such that each element belongs to $\femspacecurve{i}$.
Let $\{\basiscurve{i}{j}\}_{j=1}^{N_i}$ be the standard basis of $\femspacecurve{i}$ for $1\leq i\leq 3$, namely $\basiscurve{i}{j}{(q_{i,k})} = {\delta_{jk}}$ holds.
We set $V(\curve{m}{}) := \bigotimes_{i=1}^3\femspacecurve{i}$ and
\begin{equation*}
    \underline{V}(\curve{m}{}) := \left\{(\X{m}{1},\X{m}{2},\X{m}{3})\in\bigotimes_{i=1}^3 \femspacecurvevector{i}\mid \X{m}{1} = \X{m}{2} = \X{m}{3}\ \mbox{on}\ \partial I\right\}.
\end{equation*}
Let $\edge{m}{i}{j}$ denote the edge $[\vertex{m}{i}{j-1},\vertex{m}{i}{j}] := \{(1-s)\vertex{m}{i}{j-1} + s\vertex{m}{i}{j}\mid 0\leq s\leq 1\}$.
We define the normal vector to each edge of $\curve{m}{i}$ by 
\begin{equation*}
    \normal{m}{i}{j-\frac{1}{2}} := \frac{(\vertex{m}{i}{j} - \vertex{m}{i}{j-1})^\perp}{|\edge{m}{i}{j}|}\ \ \mbox{for}\ \ 1\leq j\leq N_i,
\end{equation*}
where $\binom{a}{b}^\perp = \binom{-b}{a}$ denotes the anti-clockwise rotation of $\vec{p}$ through $\frac\pi2$ and $|\edge{m}{i}{j}|$ is the length of the interval $\edge{m}{i}{j}$.
Let {$\normalnoindex{m}{i}$} be the normal vector field on $\curve{m}{i}$ which is equal to $\normal{m}{i}{j-\frac{1}{2}}$ on each edge $\edge{m}{i}{j}\ (1\leq j\leq N_i)$.

For two piecewise continuous functions on {$I$}, which may jump across the {points} ${q_{i,j}}\ (1\leq j\leq N_i)$, we define the mass lumped inner
product 
\begin{equation*}
    \left<u,v\right>^h_{\curve{m}{i}} := \frac{1}{2}\sum_{j = 1}^{N_i}|\edge{m}{i}{j}|((u\, v)({q^-_{i,j}}) + (u\, v)({q^+_{i,j-1}})),
\end{equation*}
where $u({q^{-}_{i,j}}) := \lim_{{[q_{i,j-1},q_{i,j}]}\ni y\to{q_{i,j}}}u(y)$ and $u({q^{+}_{i,j}}) := \lim_{{[q_{i,j},q_{i,j+1}]}\ni y\to{q_{i,j}}}u(y)$ for each $1\leq i\leq 3$.
We extend these definitions to vector- and tensor-valued functions.
Moreover, we define
    \begin{equation*}
        \innerproduct{u}{v}{(h)}{\curve{m}{}} := \sum_{i=1}^3 \innerproduct{u\!\mid_{\curve{m}{i}}}{v\!\mid_{\curve{m}{i}}}{(h)}{\curve{m}{i}},
    \end{equation*}
where here and throughout, 
the notation $\cdot^{(h)}$ means an expression with or
without the superscript $h$.
The vertex normals $\weightnormalnoindex{m}{i}\in \femspacecurvevector{i}$ {on $\curve{m}{i}$} are defined through the lumped $L^2$ projection
\begin{equation}\label{eq:weightednormaldef}
    \innerproduct{\weightnormalnoindex{m}{i}}{\vec{\xi}}{h}{\curve{m}{i}} = \innerproduct{\normalnoindex{m}{i}}{\vec{\xi}}{}{\curve{m}{i}}\ \ \mbox{for}\ \ \forall{\vec{\xi}}\in {\femspacecurvevector{i}},
\end{equation}
see \cite[Definition~51]{BarrettGarckeRobertBook2020}.

We make an assumption on the discrete vertex normals $\weightnormalnoindex{m}{}$, following \cite[Assumption~$\mathcal{A}$]{BarrettGarckeRobert2007} and \cite[Assumption~108]{BarrettGarckeRobertBook2020}, which will guarantee well-posedness of the system of linear equations:
\begin{asm} \label{asm:1}
Assume that 
$\operatorname{span}{\{\weightnormal{m}{i}{j}\}_{1\leq j\leq N_i-1}} \not= 
\{\vec0\}$ for $1 \leq i \leq 3$ and
    \begin{equation*} 
        \operatorname{span}\left\{
\sum_{(i,j,k)\in\Lambda}\innerproduct{\weightnormalnoindex{m}{i}}{\varphi_j - \varphi_k}{(h)}{\curve{m}{i}}
 \mid\ \bv \varphi\in\femspacebulkvectorsum\right\} = \bR^2.
    \end{equation*}
Here, for $\vec{\xi}\in \femspacecurvevector{i}$ and $\varphi \in S^m$,
we use the slight abuses of notation
$\innerproduct{\vec\xi}{\varphi}{}{\curve{m}{i}}
= \int_{\curve{m}{i}} \vec\xi \varphi \dH{1}$ and
$\innerproduct{\vec\xi}{\varphi}{h}{\curve{m}{i}}
= \int_{\curve{m}{i}} \pi^h[\vec\xi \varphi] \dH{1}$.
\end{asm}
We remark that the first condition basically means that each of the three
curves has at least one nonzero inner vertex normal, something that can only be
violated in very pathological cases. The proof of
Theorem~\ref{thm:exist_unique} shows that it is actually sufficient to require
this for just two out of the three curves, but for simplicity we prefer to
state the stronger assumption.
The second condition in Assumption~\ref{asm:1}, on the other hand,
is a very mild constraint on the interaction between bulk and interface meshes.
In fact, it can only be violated if all the vectors in the set are linearly
dependent, which happens, for example, if the three curves are straight lines 
that lie on top of each other.

Given $(\X{m}{},\kappa^m_\sigma)\in \underline{V}(\curve{m}{})\times V(\curve{m}{})$, 
we find $(\chemicalDiscVec{m+1}{}, \X{m+1}{}, \kappa^{m+1}_\sigma)\in \femspacebulkvectorsum\times \underline{V}(\curve{m}{})\times V(\curve{m}{})$ such that the following conditions hold:
\newline\newline\noindent
\textbf{[Motion law]} For all $\bv\varphi\in\femspacebulkvectorsum$,
\begin{subequations} \label{eq:FEA}
\begin{equation}\label{eq:FE-ML-linear}
\innerproduct{\nabla \chemicalDiscVec{m+1}{}}{\nabla\bv\varphi}{}{\Omega} - 
\sum_{(i,j,k)\in\Lambda}\innerproduct{\Glinear{m}{i}}{\varphi_j - \varphi_k}{(h)}{\curve{m}{i}}
= 0.
\end{equation}
\textbf{[Gibbs--Thomson law]} For all $\xi\in V(\curve{m}{})$,
\begin{equation}\label{eq:FE-GTL-linear}
\innerproduct{\kappa^{m+1}_\sigma}{\xi}{h}{\curve{m}{}} -
\sum_{(i,j,k)\in\Lambda}
\innerproduct{\chemicalDisc{m+1}{j} - \chemicalDisc{m+1}{k}}{\xi_i}{(h)}{\curve{m}{i}} 
= 0.
\end{equation}
\textbf{[Curvature vector]} For all $\vec\eta\in\underline{V}(\curve{m}{})$,
\begin{equation}\label{eq:FE-CV-linear}
    \innerproduct{\kappa^{m+1}_\sigma\weightnormalnoindex{m}{}}{\vec{\eta}}{h}{\curve{m}{}} + \innerproduct{\sigma\sgrad\X{m+1}{}}{\sgrad\vec{\eta}}{}{\curve{m}{}} = 0.
\end{equation}
\end{subequations}

We stress that \eqref{eq:FEA} encodes two different schemes:
One that uses mass lumping in the two bulk-surface terms in 
\eqref{eq:FE-ML-linear} and \eqref{eq:FE-GTL-linear}, and one
that uses true integration in both. The interpolation operator $\pi^h_i$
in \eqref{eq:FE-ML-linear} is superfluous in the former case, but necessary
for the stability proof of the latter. Writing \eqref{eq:FEA} as above allows
for a compact presentation. Observe that the implementation of the
scheme with mass-lumping is far easier, since there bulk finite element
functions only need to be evaluated at the vertices of the curve network.
We refer to \cite{BarrettGarckeRobert2010} for more details.

\begin{thm}[Existence and uniqueness]\label{thm:exist_unique}
Let Assumption~\ref{asm:1} hold and let $m\geq0$. 
Then there exists a unique solution
$(\chemicalDiscVec{m+1}{}, \X{m+1}{}, \kappa^{m+1}_\sigma)\in
\femspacebulkvectorsum\times \underline{V}(\curve{m}{})\times V(\curve{m}{})$
to \eqref{eq:FEA}.
\end{thm}
\begin{proof}
Since \eqref{eq:FE-ML-linear}, \eqref{eq:FE-GTL-linear} and 
\eqref{eq:FE-CV-linear} is a
linear system with the same number of unknowns and equations, existence follows
from uniqueness. To show the latter, it is sufficient to prove that only the
zero solution solves the homogeneous system. Hence let
$(\chemicalDiscVec{}{}, \X{}{}, \kappa_\sigma)\in \femspacebulkvectorsum\times \underline{V}(\curve{m}{})\times V(\curve{m}{})$ be such that 
\begin{subequations} \label{eq:FEAproof}
    \begin{align}
& \tau_m \innerproduct{\nabla \chemicalDiscVec{}{}}{\nabla\bv\varphi}{}{\Omega} - \sum_{(i,j,k)\in\Lambda}\innerproduct{\pi^h_i\left[\X{}{i}\cdot\weightnormalnoindex{m}{i}\right]}{\varphi_j - \varphi_k}{(h)}{\curve{m}{i}} = 0\ \ \forall\bv \varphi\in\femspacebulkvectorsum, \label{eq:FE-ML} \\
& \innerproduct{ \kappa_{\sigma}}{\xi}{h}{{\curve{m}{}}} - \sum_{(i,j,k)\in\Lambda}\innerproduct{\chemicalDisc{}{j} - \chemicalDisc{}{k}}{\xi_{i}}{(h)}{\curve{m}{i}}
 = 0\ \ \forall\xi\in {\femspacegamma}, \label{eq:FE-GTL} \\ &
        \innerproduct{\kappa_{\sigma}\weightnormalnoindex{m}{}}{\vec{\eta}}{h}{\curve{m}{}} + \innerproduct{\sigma\sgrad\X{}{}}{\sgrad\vec{\eta}_{}}{}{\curve{m}{}} = 0\ \ \forall\vec{\eta}\in {\underline{V}(\curve{m}{})}.
\label{eq:FE-CV}
\end{align}
\end{subequations}
    Choosing $\bv \varphi = \chemicalDiscVec{}{}$ in \eqref{eq:FE-ML}, $\xi = \pi^h[\X{}{}\cdot\weightnormalnoindex{m}{}]$ in \eqref{eq:FE-GTL} and
$\vec{\eta} = \X{}{} \in \underline{V}(\curve{m}{})$ in $\eqref{eq:FE-CV}$ 
gives
    \begin{equation}\label{eq:U_4}
0= \tau_m\|\nabla \chemicalDiscVec{}{}\|^2_2 - \innerproduct{\kappa_{\sigma}\weightnormalnoindex{m}{}}{\X{}{}}{h}{\curve{m}{}} = \tau_m\|\nabla \chemicalDiscVec{}{}\|^2_2 +\innerproduct{\sigma\sgrad\X{}{}}{\sgrad\X{}{}}{}{\curve{m}{}}.
    \end{equation}
    Thus, we see that $\chemicalDisc{}{i} = C_i\in\bR$ and $\X{}{i} = \X{c}{}\in\bR^2$, for $1\leq i\leq 3$, are constant functions, with $\sum_{i=1}^3 C_i=0$.
We deduce from \eqref{eq:FE-ML} that
\begin{equation} \label{eq:Xc_sum}
\X{c}{} \cdot \sum_{(i,j,k)\in\Lambda}\innerproduct{\weightnormalnoindex{m}{i}}{\varphi_j - \varphi_k}{(h)}{\curve{m}{i}} = 0 \quad
\forall\bv \varphi\in\femspacebulkvectorsum,
\end{equation}
and so the second condition in Assumption \ref{asm:1} yields that 
$\X{}{} = \X{c}{} = \vec{0}$.
Moreover, it follows from $\eqref{eq:FE-GTL}$ that $\kappa_{\sigma,1} = C_3 - C_2, \kappa_{\sigma,2} = C_1 - C_3$, and $\kappa_{\sigma,3} = C_2 - C_1$ are also equal to constants.
We now choose as a test function in \eqref{eq:FE-CV} the function
$\vec\eta \in\underline{V}(\curve{m}{})$ with
$\vec\eta_i =\kappa_{{\sigma,i}} \vec z^m_i$, where
$\vec z^m_i = \sum_{j=1}^{N_i-1} \vec \omega^m_{i,j} \Phi_{i,j}$ for $i=1,2,3$.
Hence, we obtain 
\[
\sum_{i=1}^3 (\kappa_{\sigma,i})^2 
\innerproduct{\vec z^m_i}{\vec z^m_i}{h}{\curve{m}{i}} = 0,
\]
so that the first condition in Assumption~\ref{asm:1} implies that 
$\kappa_{\sigma,i} = 0$ for $i=1,2,3$, i.e.\ $\kappa_\sigma = 0$.
Since $C_1=C_2=C_3$, they must all be zero, and so
$\bv W = \bvzero$ also follows.
Hence{,} we have shown the existence of a unique solution
$(\chemicalDiscVec{m+1}{}, \X{m+1}{}, \kappa^{m+1}_\sigma)\in
\femspacebulkvectorsum\times \underline{V}(\curve{m}{})\times V(\curve{m}{})$
to \eqref{eq:FEA}.
\end{proof}

Before proving the stability of our scheme, we recall the following lemma from \cite[Lemma~57]{BarrettGarckeRobertBook2020} without the proof.
\begin{lem}\label{lem:HandbookLemma57}
    Let $\curve{h}{}$ be polygonal curve in $\mathbb{R}^2$. Then, for any $\X{}{}\in \underline{V}(\curve{h}{})$, it holds that
    \begin{equation*}
        \innerproduct{\sgrad\X{}{}}{\sgrad(\X{}{} - \vec{Id})}{}{\curve{h}{}} \geq |\X{}{}(\curve{h}{})| - |\curve{h}{}| + \innerproduct{|\sgrad\X{}{}| - 1}{|\sgrad\X{}{}| - 1}{}{\curve{h}{}}{,}
    \end{equation*}
    {where $|\curve{h}{}|$ and $|\X{}{}(\curve{h}{})|$ are the lengths of $\curve{h}{}$ and $\X{}{}(\curve{h}{})$, respectively}.
\end{lem}
\begin{thm}[Unconditional stability]\label{thm:StabilityLinear}
Let $m\geq0$ and let
$(\chemicalDiscVec{m+1}{},\X{m+1}{},\kappa^{m+1}_\sigma)\in \femspacebulkvectorsum\times\underline{V}(\curve{m}{})\times V(\curve{m}{})$ be a solution to \eqref{eq:FEA}. Then the following estimate is satisfied:
    \begin{equation*}
        |\curve{m+1}{}|_\sigma + \tau_m\|\nabla \chemicalDiscVec{m+1}{}\|^2_2 \leq |\curve{m}{}|_\sigma,
    \end{equation*}
where we recall that 
$|\curve{m}{}|_\sigma = \sum_{i=1}^3\sigma_i|\curve{m}{i}|$.
\end{thm}
\begin{proof}
    We choose $\bv \varphi = \chemicalDiscVec{m+1}{}$ in $\eqref{eq:FE-ML-linear}$ 
    and $\xi = \pi^h\left[\frac{\X{m+1}{}-\X{m}{}}{\tau_m}\cdot\weightnormalnoindex{m}{}\right]$ in $\eqref{eq:FE-GTL-linear}$ to obtain
    \begin{align}
        \|\nabla\chemicalDiscVec{m+1}{}\|^2_2 
        & = \innerproduct{\kappa^{m+1}_{\sigma}}{\Glinear{m}{}}{h}{\curve{m}{}}
        = \frac1{\tau_m} \innerproduct{\kappa_{\sigma}^{m+1}\vec\omega^m}{\X{m+1}{} -
        \X{m}{}}{h}{\curve{m}{}}. \label{eq:StabilityLinear3}
    \end{align}
    Choosing $\vec{\eta} = \X{m+1}{} - \X{m}{}\in\underline{V}(\curve{m}{})$ in $\eqref{eq:FE-CV-linear}$, gives
    \begin{equation}\label{eq:StabilityLinear4}
        \innerproduct{\kappa_{\sigma}^{m+1}\vec\omega^m}{\X{m+1}{} - \X{m}{}}{h}{\curve{m}{}} = -\innerproduct{\sigma\sgrad\X{m+1}{}}{\sgrad(\X{m+1}{} - \X{m}{})}{}{\curve{m}{}}.
    \end{equation}
    We compute from \eqref{eq:StabilityLinear3} and \eqref{eq:StabilityLinear4}, on recalling Lemma~\ref{lem:HandbookLemma57}, that
    \[
        \tau_m\|\nabla \chemicalDiscVec{m+1}{}\|^2_2 =
-\innerproduct{\sigma\sgrad\X{m+1}{}}{\sgrad(\X{m+1}{} - \X{m}{})}{}{\curve{m}{}}
\leq |\curve{m}{}|_\sigma - |\curve{m+1}{}|_\sigma.
    \]
This proves the desired result.
\end{proof}
\begin{rem}
Similar stability results for the discretization of 
other gradient flows of curve networks can
be found in \cite{BarrettGarckeNuernbergSurfaceDiffusion2007,BarrettGarckeRobert2007}. For the two-phase Mullins--Sekerka problem, related stability results
were derived in \cite{BarrettGarckeRobert2010,Nurnberg202203}.
\end{rem}

\section{Solution of the linear system}\label{sec:MatrixForm}
In this section, we discuss solution methods for the systems of
linear equations arising from \eqref{eq:FEA} at each time level.
To this end, we make use of ideas from 
\cite{BarrettGarckeRobert2007,BarrettGarckeRobert2010}.
Here the crucial idea is to avoid having to work with the 
trial and test spaces $\underline{V}(\curve{m}{})$ and
$\femspacebulkvectorsum$ directly, and rather employ
a technique that is similar to a standard treatment of 
periodic boundary conditions for ODEs and PDEs.
In particular, following \cite[(2.44)]{BarrettGarckeRobert2007} we
introduce the orthogonal projections 
$\mathcal{P}:[V(\curve{m}{})]^2\to \underline{V}(\curve{m}{})$ 
and $\mathcal{Q}:\femspacebulkvector\to \femspacebulkvectorsum$.
On letting $\bvone = (1,1,1)^T$, it is easy to see that for 
$\bv W \in\femspacebulkvector$ it holds that
    \begin{equation*}
        \mathcal{Q}\chemicalDiscVec{}{} = \chemicalDiscVec{}{} -
\frac{\chemicalDiscVec{}{}\cdot\bvone}{\bvone\cdot\bvone}\, \bvone
    \end{equation*}
point-wise in $\overline\Omega$.

Now, given $\X{m}{} \in \underline{V}(\curve{m}{})$, let $(\chemicalDiscVec{m+1}{},\kappa^{m+1}_\sigma,\X{m}{}+\delta\X{m+1}{})\in\femspacebulkvectorsum\times V(\curve{m}{})\times \underline{V}(\curve{m}{})$ be the unique solution to \eqref{eq:FEA} whose existence has been proven in Theorem~\ref{thm:exist_unique}. 
Let $N := \sum_{i=1}^3(N_i+1)$ be the sum of the vertices on each individual
curve, 
and let 
$\spacemeshpointsnum$ be the number of vertices of the mesh $\mathfrak{T}^m$ inside $\overline\Omega$. 
{From} now on, as no confusion can arise, we identify
$(\chemicalDiscVec{m+1}{},\kappa^{m+1}_\sigma,\delta\X{m+1}{})$ with their
vectors of coefficients with respect to the bases $\{\basisbulk{m}{i}\}_{1\leq i\leq\spacemeshpointsnum}$ and $\{\{\basiscurve{i}{j}\}_{1\leq j\leq N_i}\}_{i=1}^3$ of the unconstrained spaces $\femspacebulkvector$ and 
$V(\curve{m}{})$. In addition, we let $\metric{P} : (\bR^2)^N \to
\{(\vec{z}_1,\vec{z}_2,\vec{z}_3) \in (\bR^2)^N :
[\vec{z}_1]_0 = [\vec{z}_2]_0 = [\vec{z}_3]_0,
[\vec{z}_1]_{N_1}=[\vec{z}_2]_{N_2}=[\vec{z}_3]_{N_3}\}$ be the Euclidean
space equivalent of $\mathcal{P}$, and similarly for the equivalent
$Q : (\mathbb{R}^3)^{\spacemeshpointsnum}\to \mathbb{W} := \{(v_1,v_2,v_3)^T\in(\mathbb{R}^3)^{\spacemeshpointsnum}\ \mid\ v_1+v_2+v_3 = 0\in\mathbb{R}^{\spacemeshpointsnum}\}$ of $\mathcal{Q}$.

Then the solution to \eqref{eq:FEA} can be written as
$(Q \chemicalDiscVec{m+1}{}, \kappa^{m+1}_\sigma,
 \X{m}{} + \metric{P}\delta\X{m+1}{})$ for any solution of the
linear system 
\begin{equation}\label{eq:MatrixForm}
    \begin{pmatrix}
        \tau_m QA_\Omega Q & O & Q\vec{N}_{\Omega,\Gamma}^T\metric{P} \\
        B_{\Omega,\Gamma}Q & C_\Gamma & O \\
        O & \metric{P}\vec{D}_\Gamma & \metric{P}\metric{E_\Gamma}\metric{P}
    \end{pmatrix}
    \begin{pmatrix}
        \chemicalDiscVec{m+1}{}\\
        \kappa^{m+1}_\sigma\\
        \delta\X{m+1}{}
    \end{pmatrix} = 
    \begin{pmatrix}
        O\\
        O\\
        -\metric{P}\metric{E_\Gamma}\metric{P}\X{m}{}
    \end{pmatrix}{,}
\end{equation}
where $A_\Omega\in\mathbb{R}^{3\spacemeshpointsnum\times 3\spacemeshpointsnum},\vec{N}_{\Omega,\Gamma}\in(\mathbb{R}^2)^{N\times 3\spacemeshpointsnum},B_{\Omega,\Gamma}\in\mathbb{R}^{N\times 3\spacemeshpointsnum},C_\Gamma\in\mathbb{R}^{N\times N},\vec{D}_\Gamma\in(\mathbb{R}^2)^{N\times N}$ and $\metric{{E}_\Gamma}\in(\mathbb{R}^{2\times 2})^{N\times N}$ are defined by
\begin{multline*}
    A_\Omega := 
    \begin{pmatrix}
        A & O & O  \\
        O & A & O \\
        O & O & A
    \end{pmatrix},\ 
   \vec{N}_{\Omega,\Gamma} :=
    \begin{pmatrix}
        O & \vec{N}_1 & -\vec{N}_1 \\
        -\vec{N}_2 & O & \vec{N}_2 \\
        \vec{N}_3 & -\vec{N}_3 & O
    \end{pmatrix},\ 
    B_{\Omega,\Gamma} := 
    \begin{pmatrix}
        O & B_1 & -B_1 \\
        -B_2 & O & B_2 \\
        B_3 & -B_3 & O \\
    \end{pmatrix},\\ 
    C_\Gamma := 
    \begin{pmatrix}
        C_1 & O & O \\
        O & C_2 & O \\
        O & O & C_3 \\
    \end{pmatrix},\ 
    \vec{D}_\Gamma :=
    \begin{pmatrix}
        \vec{D}_1 & O & O \\
        O & \vec{D}_2 & O \\
        O & O & \vec{D}_3
    \end{pmatrix},\ 
    \metric{{E}_\Gamma} :=
    \begin{pmatrix}
        \sigma_1\metric{E_1} & O & O \\
        O & \sigma_2\metric{E_2} & O \\
        O & O & \sigma_3\metric{E_3} \\
    \end{pmatrix}{,}
\end{multline*}
with
\begin{equation}\label{eq:submatrices}
    \begin{array}{ll}
        \left[A\right]_{i,j} := \left<\nabla\Psi^m_j,\nabla\Psi^m_i\right>,  & \left[\vec{N}_c\right]_{l,i} := \innerproduct{\Phi^m_{c,l}}{\Psi^m_i}{{(h)}}{\curve{m}{c}}\weightnormal{m}{c}{l},\\
        \left[B_c\right]_{k,j} := \innerproduct{\Psi^m_j}{\Phi^m_{c,k}}{{(h)}}{\curve{m}{c}}, & \left[C_c\right]_{k,l} := \innerproduct{\Phi^m_{c,l}}{\Phi^m_{c,k}}{h}{\curve{m}{c}}, \\
        \left[\vec{D}_c\right]_{k,l} := \innerproduct{\Phi^m_{c,l}}{\Phi^m_{c,k}}{h}{\curve{m}{c}}\weightnormal{m}{c}{l}, & \left[\metric{{E}_c}\right]_{k,l} := \innerproduct{\sgrad\Phi^m_{c,l}}{\sgrad\Phi^m_{c,k}}{}{\curve{m}{c}},\\
    \end{array}
\end{equation}
for each $1\leq c\leq 3$.
The advantage of the system \eqref{eq:MatrixForm} over a naive implementation
of \eqref{eq:FEA} is that complications due to nonstandard finite element
spaces are completely avoided. A disadvantage is, however, that the system
\eqref{eq:MatrixForm} is highly singular, in that due to the presence 
of the projections the dimension of its kernel is larger than the dimension of
the scalar bulk finite element space $S^m$.
This
makes it difficult to solve \eqref{eq:MatrixForm} in practice. A more practical
formulation can be obtained by eliminating one of the components of
$\chemicalDiscVec{m+1}{}$ completely. 
In particular, on recalling that $\chemicalDiscVec{m+1}{} \cdot \bvone = 0$, we can reduce the unknown variables $\chemicalDiscVec{m+1}{}\in(\mathbb{R}^3)^{\spacemeshpointsnum}$ to $(W^{m+1}_1,W^{m+1}_2)\in(\mathbb{R}^2)^{\spacemeshpointsnum}$ 
    by introducing the linear map $\widehat{Q}:(\mathbb{R}^2)^{\spacemeshpointsnum}\to \mathbb{W}\subset(\mathbb{R}^3)^{\spacemeshpointsnum}$ defined by
    \begin{equation*}
        \widehat{Q} := \begin{pmatrix}
            I_{\spacemeshpointsnum} & O \\
            O & I_{\spacemeshpointsnum}\\
            -I_{\spacemeshpointsnum} & -I_{\spacemeshpointsnum}\\
        \end{pmatrix},
    \end{equation*}
    where $I_M$ denotes the identity matrix of size $M$ for $M\in\mathbb{N}$. 

Then the solution to \eqref{eq:FEA} can be written as
$(\widehat Q \widehat{\bv W}^{m+1}, \kappa^{m+1}_\sigma,
 \X{m}{} + \metric{P}\delta\X{m+1}{})$ for any solution of the reduced
linear system 
    \begin{equation}\label{eq:MatrixFormReduced}
        \begin{pmatrix}
            \widehat{A}_\Omega & O & \widehat{N}^T_{\Omega,\Gamma}\metric{P} \\
            \widehat{B}_{\Omega,\Gamma} & C_\Gamma & O \\
            O & \metric{P}\vec{D}_\Gamma & \metric{P}\metric{{E}_\Gamma}\metric{P}
        \end{pmatrix}
        \begin{pmatrix}
            \widehat{\bv W}^{m+1}\\
            \kappa_\sigma^{m+1}\\
            \delta\X{m+1}{}
        \end{pmatrix} = 
        \begin{pmatrix}
            O \\
            O \\
            -\metric{P}\metric{{E}_\Gamma}\metric{P}\X{m}{}
        \end{pmatrix},
    \end{equation}
    where
    \begin{equation*}
        \widehat{A}_\Omega :=
        \begin{pmatrix}
            A & O \\
            O & A \\
        \end{pmatrix} ,\ 
        \widehat{B}_{\Omega,\Gamma} := B_{\Omega,\Gamma}\widehat{Q} = 
        \begin{pmatrix}
            B_1 & 2B_1 \\
            -2B_2 & -B_2 \\
            B_3 & -B_3
        \end{pmatrix},\ 
        \widehat{N}_{\Omega,\Gamma} :=
        \begin{pmatrix}
            O & \vec{N}_1 \\
            -\vec{N}_2 & O \\
            \vec{N}_3 & -\vec{N}_3
        \end{pmatrix}.
    \end{equation*}
In contrast to \eqref{eq:MatrixForm}, the kernel of 
\eqref{eq:MatrixFormReduced} is small. In fact, it has dimension $8$
due to the fact that $\metric{P}$ has a kernel of dimension $8$.
Hence, iterative solution methods, combined with good
preconditioners, work very well to solve \eqref{eq:MatrixFormReduced} in
practice.
For our numerical results in Section~\ref{sec:NumericalResults}, 
below, we employ
a GMRes iterative solver with least squares solution of the block matrix in 
\eqref{eq:MatrixFormReduced} without $\metric{P}$ 
as preconditioners.

\section{Obtaining a fully discrete area conservation property}\label{sec:ImprovementOfTheScheme}

Although the linear scheme \eqref{eq:FEA} introduced in
Section~\ref{sec:FEM_approx} can be shown to be unconditionally
stable, recall Theorem~\ref{thm:StabilityLinear}, in general the areas occupied by the
discrete approximations of the three phases will not be conserved.
In this section, we state how to modify the previously introduced scheme 
\eqref{eq:FEA} in such a way, that it satisfies both of the structure defining 
properties from Section~\ref{sec:MathematicalProperties}. To this end,
we follow the discussion in \cite[\S3]{Nurnberg202203} in order to obtain
an exact area preservation property on the fully discrete level.
We remark that our approach hinges on ideas first presented for area-conserving
geometric flows for closed curves in \cite{JiangL21,BaoZ21}.
See also \cite[\S3.2]{BaoGarckeNuernbergZhao2022} for related work in
the context of the surface diffusion flow for curve networks with 
triple junctions.

Let us define families of polygonal curves $\{\curve{h}{i}(t)\}_{t\geq 0}$,
$i=1,2,3$, that {are} parameterized by the time variable. 
In particular, for each $0\leq m\leq M$, $t\in[t_{m},t_{m+1}]$ and $1\leq i\leq 3$, we define the polygonal curve $\curve{h}{i}(t)$ by
\begin{equation*}
    \curve{h}{i}(t) := \frac{t_{m+1} - t}{\tau_m}\curve{m}{i} + \frac{t - t_m}{\tau_m}\curve{m+1}{i}.
\end{equation*}
Precisely speaking, the vertices of $\curve{h}{i}(t)$ {are} defined as follows:
\begin{equation*}
    \vertex{h}{i}{j}(t) := \frac{t_{m+1}-t}{\tau_m}\vertex{m}{i}{j} + \frac{t - t_m}{\tau_m}\vertex{m+1}{i}{j}\ \ \mbox{for}\ \ 0\leq j\leq N_i,
\end{equation*}
while we write each edge of $\curve{h}{i}(t)$ as $\edge{h}{i}{j}(t) := [\vertex{h}{i}{j-1}(t),\vertex{h}{i}{j}(t)]$ {for $1\leq i\leq 3$ and} $1\leq j\leq N_i$.

\begin{lem}\label{lem:AreaDifferenceDiscrete}
For each $m\geq0$ and $(i,j,k)\in\Lambda$, it holds that
    \begin{equation*}
    |\Omega^{m+1}_i| - |\Omega^m_i| = \innerproduct{\X{m+1}{j} - \X{m}{j}}{\vertexnormalnoindex{m+\frac{1}{2}}{j}}{}{\curve{m}{j}} - \innerproduct{\X{m+1}{k}-\X{m}{k}}{\vertexnormalnoindex{m+\frac{1}{2}}{k}}{}{\curve{m}{k}},
    \end{equation*}
    where
    \begin{equation}\label{eq:normalInterporated}
        \vertexnormalnoindex{m+\frac{1}{2}}{i} := \frac{1}{\tau_m|\edge{m}{i}{j}|}\int_{t_m}^{t_{m+1}}\normalnoindex{h}{}|\edge{h}{i}{j}(t)|\dt\ \ \mbox{on}\ \ \sigma^m_{i,j}\ \ \mbox{for}\ \ 1\leq j\leq N_i,
    \end{equation}
\end{lem}
\begin{proof}
The desired result for $i=1,2$ is shown in
\cite[Lemma~3.1]{BaoGarckeNuernbergZhao2022}, and the result for
$i=3$ follows analogously on noting that $\partial\Omega$ is fixed.
\end{proof}

Now the weighted vertex normal vector $\weightnormalnoindex{m+\frac{1}{2}}{i}\in\femspacecurvevector{i}$, for $i=1,2,3$, associated with $\vertexnormalnoindex{m+\frac{1}{2}}{i}$ is defined through the following formula:
\begin{equation}\label{eq:WeightedNormalAndEdgeNormal}
    \innerproduct{\weightnormalnoindex{m+\frac{1}{2}}{i}}{\vec{\xi}}{h}{\curve{m}{i}} = \innerproduct{\vertexnormalnoindex{m+\frac{1}{2}}{i}}{\vec{\xi}}{}{\curve{m}{i}}\ \ \forall\vec{\xi}\in {\femspacecurvevector{i}}.
\end{equation}
{
    Consequently, we obtain a nonlinear system with the aid of $\weightnormalnoindex{m+\frac{1}{2}}{i}$:
Given $(\X{m}{},\kappa^m_\sigma)\in \underline{V}(\curve{m}{})\times V(\curve{m}{})$, 
find $(\chemicalDiscVec{m+1}{}, \X{m+1}{}, \kappa^{m+1}_\sigma)\in \femspacebulkvectorsum\times \underline{V}(\curve{m}{})\times V(\curve{m}{})$ such that
\newline\newline
}
\textbf{[Motion law]} For all $\bv \varphi\in\femspacebulkvectorsum$,
\begin{subequations}\label{eq:FEAnl}
\begin{equation}\label{eq:FE-ML-dash}
   \innerproduct{\nabla \chemicalDiscVec{m+1}{}}{\nabla\bv \varphi}{}{\Omega} - \sum_{(i,j,k)\in\Lambda}\innerproduct{\pi^h_i\left[\frac{\X{m+1}{i} - \X{m}{i}}{\tau_m}\cdot\weightnormalnoindex{m+\frac{1}{2}}{i}\right]}{\varphi_j - \varphi_k}{(h)}{\curve{m}{i}} = 0.
\end{equation}
\textbf{[Gibbs--Thomson law]} For all $\xi\in V(\curve{m}{})$,
\begin{equation}\label{eq:FE-GTL-dash}
\innerproduct{\kappa^{m+1}_\sigma}{\xi}{h}{\curve{m}{}} - \sum_{(i,j,k)\in\Lambda}\innerproduct{\chemicalDisc{m+1}{j} - \chemicalDisc{m+1}{k}}{\xi_i}{(h)}{\curve{m}{i}} = 0.
\end{equation}
\textbf{[Curvature vector]} For all {$\vec{\eta}\in \underline{V}(\curve{m}{})$},
\begin{equation}\label{eq:FE-CV-dash}
    {\innerproduct{\kappa^{m+1}_\sigma\weightnormalnoindex{m+\frac{1}{2}}{}}{\vec{\eta}}{h}{\curve{m}{}} + \innerproduct{\sigma\sgrad\X{m+1}{}}{\sgrad\vec{\eta}}{}{\curve{m}{}}} = 0.
\end{equation}
\end{subequations}

We can now prove the area preserving property of each domain surrounded by the polygonal curve in the discrete level.

\begin{thm}[Area preserving property for the discrete scheme]\label{thm:AP_nonlinear}
    Let $m\geq0$ and let \linebreak
$(\chemicalDiscVec{m+1}{},\X{m+1}{},\kappa^{m+1}_\sigma)\in\femspacebulkvectorsum\times\underline{V}(\curve{m}{})\times V(\curve{m}{})$ be a solution of $\eqref{eq:FEAnl}$.
    Then, for each $1\leq j\leq 3$ it holds that
    \begin{equation*}
        |\Omega^{m+1}_j| = |\Omega^m_j|.
    \end{equation*}
\end{thm}
\begin{proof}
We will argue similarly to the proof of \cite[Theorem~3.3]{Nurnberg202203}.
Choosing $\bv \varphi = (-\frac23,\frac{1}{3},\frac{1}{3})^T\in\femspacebulkvectorsum$ in $\eqref{eq:FE-ML-dash}$, we obtain from
\eqref{eq:WeightedNormalAndEdgeNormal} and 
Lemma~\ref{lem:AreaDifferenceDiscrete} that
    \begin{align*}
        0 & = \innerproduct{\pi^{h}_{{3}}\left[\frac{\X{m+1}{3} - \X{m}{3}}{\tau_m}\cdot\weightnormalnoindex{m+\frac{1}{2}}{3}\right]}{1}{{(h)}}{\curve{m}{3}} + \innerproduct{\pi^{h}_{{2}}\left[-\frac{\X{m+1}{2}-\X{m}{2}}{\tau_m}\cdot\weightnormalnoindex{m+\frac{1}{2}}{2}\right]}{1}{{(h)}}{\curve{m}{2}}\\&
        =\innerproduct{\frac{\X{m+1}{3}-\X{m}{3}}{\tau_m}}{\weightnormalnoindex{m+\frac{1}{2}}{3}}{h}{\curve{m}{3}} + \innerproduct{-\frac{\X{m+1}{2}-\X{m}{2}}{\tau_m}}{\weightnormalnoindex{m+\frac{1}{2}}{2}}{h}{\curve{m}{2}}
        \\& = \innerproduct{\frac{\X{m+1}{3}-\X{m}{3}}{\tau_m}}{\vertexnormalnoindex{m+\frac{1}{2}}{3}}{}{\curve{m}{3}} + \innerproduct{\frac{\X{m+1}{2}-\X{m}{2}}{\tau_m}}{-\vertexnormalnoindex{m+\frac{1}{2}}{2}}{}{\curve{m}{2}} = \frac{1}{\tau_m}(|\Omega^{m+1}_1| - |\Omega^m_1|).
    \end{align*}
This yields the desired result for $k=1$. The other cases can be treated
analogously.
\end{proof}
\begin{thm}[Stability for the area-conserving scheme]\label{thm:ST_nonlinear}
Let $m\geq0$ and let
$(\chemicalDiscVec{m+1}{},\X{m+1}{},\kappa^{m+1}_\sigma)\in {\femspacebulkvectorsum}\times\underline{V}(\curve{m}{})\times V(\curve{m}{})$ be a solution of \eqref{eq:FEAnl}. Then, it holds that
    \begin{equation*}
        {|\curve{m+1}{}|_\sigma} + \tau_m\|\nabla {\chemicalDiscVec{m+1}{}}\|^2_2\leq {|\curve{m}{}|_\sigma}.
    \end{equation*}
\end{thm}
\begin{proof}
    The proof is analogous to the proof of Theorem~\ref{thm:StabilityLinear} once we replace $\weightnormalnoindex{m}{}$ by $\weightnormalnoindex{m+\frac{1}{2}}{}$.
\end{proof}
\begin{rem}
We observe that as $\weightnormalnoindex{m+\frac{1}{2}}{i}$ depends on $\X{m+1}{i}$, the scheme \eqref{eq:FEAnl} 
is no longer linear. In practice the nonlinear systems of equations arising 
at each time level of \eqref{eq:FEAnl} can be solved with the aid of a 
simple lagged iteration as mentioned in \cite[\S3]{Nurnberg202203}.
In particular, given $\Gamma^m$, let $\Gamma^{m+1,0} =
\Gamma^m$. Then for $\ell \geq 0$, and until convergence, 
define $\vec\omega^{m+\frac12,\ell}$
through \eqref{eq:WeightedNormalAndEdgeNormal} 
and \eqref{eq:normalInterporated}, 
but with $\Gamma^m$ replaced by $\Gamma^{m+1,\ell}$, and find
$(\chemicalDiscVec{m+1,\ell+1}{}, \X{m+1,\ell+1}{}, \kappa^{m+1,\ell+1}_\sigma)\in \femspacebulkvectorsum\times \underline{V}(\curve{m}{})\times V(\curve{m}{})$
such that
\begin{align*}
&
\innerproduct{\nabla \chemicalDiscVec{m+1,\ell+1}{}}{\nabla\bv \varphi}{}{\Omega} - \sum_{(i,j,k)\in\Lambda}\innerproduct{\pi^h_i\left[\frac{\X{m+1,\ell+1}{i} - \X{m}{i}}{\tau_m}\cdot\weightnormalnoindex{m+\frac{1}{2},\ell}{i}\right]}{\varphi_j - \varphi_k}{(h)}{\curve{m}{i}} = 0, \\ &
\innerproduct{\kappa^{m+1,\ell+1}_\sigma}{\xi}{h}{\curve{m}{}} - \sum_{(i,j,k)\in\Lambda}\innerproduct{\chemicalDisc{m+1,\ell+1}{j} - \chemicalDisc{m+1,\ell+1}{k}}{\xi_i}{(h)}{\curve{m}{i}} = 0,\\ &
    {\innerproduct{\kappa^{m+1,\ell+1}_\sigma\weightnormalnoindex{m+\frac{1}{2},\ell}{}}{\vec{\eta}}{h}{\curve{m}{}} + \innerproduct{\sigma\sgrad\X{m+1,\ell+1}{}}{\sgrad\vec{\eta}}{}{\curve{m}{}}} = 0,
\end{align*}
for all $\bv \varphi\in\femspacebulkvectorsum$, $\xi\in V(\curve{m}{})$,
$\vec{\eta}\in \underline{V}(\curve{m}{})$.
\end{rem}

\section{Generalization to multi-component systems}\label{sec:GeneralizationToMulti}

In order to simplify the presentation, in the previous sections we
concentrated on the simple three phase situation depicted in 
Figure~\ref{fig:problem_setting}. However, it is not difficult to generalize
our introduced finite element approximations to the general multi-phase case.
We present the details in this section, following closely the description
of the general curve network used in \cite[\S2]{BaoGarckeNuernbergZhao2022}.

\subsection{Problem setting}\label{sec:ProblemSettingMulti}
For later use, we let $\naturalset{M} := \{1,\cdots,M\}$ for each $M\in\mathbb{N}$.
    First, let us introduce counters which will be used frequently later on.
    Given a curve network $\curve{}{}(t)$, $I_C\geq1$ denotes the number of curves which are included in $\curve{}{}(t)$.
    Thus, we have $\curve{}{}(t) = \cup_{i=1}^{I_C}\curve{}{i}(t)$, where each $\curve{}{i}(t)$ is either open or closed.
    Let $I_P\geq2$ be the number of phases, i.e.\ the not necessarily connected subdomains of $\Omega$ which are separated by $\curve{}{}(t)$. This means that $\Omega\backslash\curve{}{}(t) = \bigcup_{j=1}^{I_P}\Omega_{j}(t)$.
    Finally, each endpoint of an open curve included in $\curve{}{}(t)$ is part of a triple junction.
    We write the number of triple junctions as $I_T\geq0$. 

\begin{asm}[Triple junctions]\label{asm:TripleJunction}
Every curve $\curve{}{i}(t)$, $i\in\naturalset{I_C}$, must not self-intersect, and is allowed to intersect other curves only at its boundary
$\partial\curve{}{i}(t)$. 
If $\partial\curve{}{i}(t) = \emptyset$,
then $\curve{}{i}(t)$ is called a closed curve, otherwise an open curve.
For each triple junction $\triplejunction{k}(t)$, 
$k\in\naturalset{I_T}$, there exists a unique tuple 
$(\curveindex{k}{1}, \curveindex{k}{2}, \curveindex{k}{3})$
with $1\leq\curveindex{k}{1}<\curveindex{k}{2}<\curveindex{k}{3}\leq I_C$
such that
$\triplejunction{k}(t)=\partial\curve{}{\curveindex{k}{1}}(t)\cap\partial\curve{}{\curveindex{k}{2}}(t)\cap\partial\curve{}{\curveindex{k}{3}}(t)$.
Moreover, $\cup_{i=1}^{I_C} \partial\curve{}{i}(t)\subset \cup_{k=1}^{I_T} \triplejunction{k}(t)$.
\end{asm}

\begin{asm}[Phase separation]\label{asm:DomainSeparation}
The curve network $\curve{}{}(t)$ is equipped with a matrix
$\dcmap:\{-1,0,1\}^{I_P \times I_C}$ that encodes the orientations
of the phase boundaries.
In particular, each row contains nonzero entries only for the curves that 
make up the boundary of the corresponding phase, with the sign
specifying the orientation needed for the curves normal to make it point
outwards of the phase. I.e.\ for 
$p\in\naturalset{I_P}$ and $i\in\naturalset{I_C}$ we have
\begin{equation*}
    \dcmap_{pi} = \begin{cases}
        1 & \curve{}{i}(t)\subset\pOmega_p(t)\ \mbox{and}\ \normalnoindex{}{i}\ \mbox{is outward to}\ \Omega_p(t),\\
        -1 & \curve{}{i}(t)\subset\pOmega_p(t)\ \mbox{and}\ \normalnoindex{}{i}\ \mbox{is inward to}\ \Omega_p(t),\\
        0 & \curve{}{i}(t)\not\subset\pOmega_p(t).
    \end{cases}
\end{equation*}
For every $i\in\naturalset{I_C}$, 
there exists a unique pair $(p_1,p_2)$ with $1 \leq p_1\leq p_2\leq I_P$ 
such that $\dcmap_{p_1,i} = -\dcmap_{p_2,i} = 1$.
In this situation, we say that $\Omega_{p_1}(t)$ is a neighbour to $\Omega_{p_2}(t)$.
\end{asm}

Clearly, given the matrix $\dcmap$, the boundary of $\Omega_p(t)$ can be
characterized by
\begin{equation*}
    \pOmega_p(t) = \begin{cases}
        \bigcup_{\{j:\dcmap_{pj}\neq 0\}}\curve{}{j}(t) & 1 \leq p<I_P,\\
        \pOmega\cup\bigcup_{\{j:\dcmap_{pj}\neq 0\}}\curve{}{j}(t) & p=I_P,\\
    \end{cases}
\end{equation*}
where we have assumed that $\Omega_{I_P}(t)$ is the only phase with
contact to the external boundary.

\begin{rem}[Examples]
We note that for the three-phase problem shown in
Figure~\ref{fig:problem_setting}, we have
$I_C = 3$, $I_P = 3$, $I_T = 2$, 
$(\curveindex{1}{1},\curveindex{1}{2},\curveindex{1}{3}) = (\curveindex{2}{1},\curveindex{2}{2},\curveindex{2}{3}) = (1,2,3)$
and
$\dcmap = \begin{pmatrix}0 &-1&1 \\ 1 & 0 &-1 \\ -1 &1&0 \end{pmatrix}$.
Examples for more complicated networks, and their description in our
general framework, can be found in the numerical results section, see
\S\ref{sec:83}.
\end{rem}

Under these preparations, we can generalize the system \eqref{eq:StrongForm} 
to the multi-phase case. Where no confusion can arise, we use the same notation
as before, e.g.,
$\zerosumspace = \{\bv u \in\mathbb{R}^{I_P}\mid\ \sum_{j=1}^{I_P}u_j = 0\}$. 
Given an initial curve network $\curve{0}{}$, 
our aim is to find $\chemicalVec{}{}:\Omega\to\zerosumspace$
and the evolution of a curve network $\{\curve{}{}(t)\}_{t\geq 0}$ which satisfy
\begin{equation}\label{eq:StrongForm_Multi}
    \begin{cases}
        \Delta{\bv w} = \zerovec\ \ \mbox{in}\ \ \Omega\backslash\curve{}{}(t),\ t\in[0,T],\\
        {\bv w}\cdot[\phasecharacter] = \sigma \varkappa\ \ \mbox{on}\ \ \curve{}{}(t),\ t\in[0,T],\\
        [\nabla{\bv w}]\,\vec{\nu} = -V[\phasecharacter]\ \ \mbox{on}\ \ \curve{}{}(t),\ t\in[0,T],\\
        \partial_{\normalcontainer}{\bv w} = \zerovec\ \ \mbox{on}\ \ \pOmega,\\
        \sum_{i=1}^3\sigma_{\curveindex{k}{i}}\vec{\mu}_{\curveindex{k}{i}} = \zerovec\ \ \mbox{on}\ \ \triplejunction{k}(t),\ 1\leq k\leq I_T,\ t\in[0,T],\\
        \curve{}{}(0) = \Gamma^0,
    \end{cases}
\end{equation}
where $\phasecharacter$ denotes the vector of the characteristic functions 
$\chi_j = \charfcn{\Omega_j(t)}$, $j=1,\ldots,I_P$, as before.

\begin{rem}
    The system $\eqref{eq:StrongForm_Multi}$ does not depend on the choice of normals $\normalnoindex{}{i}$, for $i\in\naturalset{I_C}$.
    Indeed, if we take $-\normalnoindex{}{i}$ as the unit normal vector to $\curve{}{i}(t)$, then the sign of $\varkappa_i$ reverses.
    On the other hand, the sign of the jump $\jump{\phasecharacter}$ is also reversed. Thus, the second law of $\eqref{eq:StrongForm_Multi}$ does not change.
    Meanwhile, the sign of the normal velocity $V_i$ in the third condition is also reversed, balancing with the sign change of $\normalnoindex{}{i}$ on the left hand side.
\end{rem}

\subsection{Weak formulation}\label{subsec:WeakFormulationMulti}
Similarly to Section~\ref{sec:WeakFormulation}, it is possible
to derive a weak formulation of the multi-phase problem 
\eqref{eq:StrongForm_Multi}.
On noting that the jump of the $j$-th characteristic function $\chi_j$
across the curve $\curve{}{i}$ is $-\dcmap_{ji}$, 
we can compute from the third condition in \eqref{eq:StrongForm_Multi} that
\begin{equation*}
    \jump{\nabla {w_j}}  \cdot \normalnoindex{}{i}= -V_i[\chi_j]
= \dcmap_{ji} V_i\ \text{ on }\ \curve{}{i}(t),
\quad 1 \leq i \leq I_C,\ 1 \leq j \leq I_P.
\end{equation*}
Hence, overall we obtain the following weak formulation:\newline\newline\noindent
\textbf{[Motion law]} For all $\bv \varphi\in \zerosumfunc\cap[H^1(\Omega)]^{I_P}$,
\begin{subequations}\label{eq:WeakFormMulti}
\begin{equation}\label{eq:WeakFormMulti_MotionLaw}
    \innerproduct{\nabla\bv w}{\nabla \bv \varphi}{}{\Omega} + \sum_{i=1}^{I_C}\sum_{j=1}^{I_P}\dcmap_{ji}\innerproduct{V_i}{\varphi_j}{}{\curve{}{i}} = 0.
\end{equation}
\textbf{[Gibbs--Thomson law]} For all $\xi\in L^2(\curve{}{})$,
\begin{equation}\label{eq:WeakFormMulti_GibbsThomsonLaw}
    \innerproduct{\varkappa_\sigma}{\xi}{}{\curve{}{}} + \sum_{i=1}^{I_C}\sum_{j=1}^{I_P}\dcmap_{ji}\innerproduct{w_j}{\xi_i}{}{\curve{}{i}} = 0.
\end{equation}
\textbf{[Curvature vector]} For all $\Vec{\eta}\in H^1(\curve{}{};\mathbb{R}^2)$ such that
$\vec{\xi}\!\mid_{\curve{}{\curveindex{k}{1}}} = \vec{\xi}\!\mid_{\curve{}{\curveindex{k}{2}}} = \vec{\xi}\!\mid_{\curve{}{\curveindex{k}{3}}}$
on $\triplejunction{k}$ for all {$k\in\naturalset{I_T}$},
\begin{equation}\label{eq:WeakFormMulti_Curvature}
    \innerproduct{\varkappa_\sigma\normalnoindex{}{}}{\Vec{\xi}}{}{\curve{}{}}  + \innerproduct{\sigma\sgrad\identity}{\sgrad\Vec{\xi}}{}{\curve{}{}} = 0.
\end{equation}
\end{subequations}

\subsection{Finite element approximations}
We now generalize our finite element approximation \eqref{eq:FEA}
to the multi-phase case. 
The necessary discrete function spaces are the obvious generalizations,
for example
\[
\femspacegammavector = \left\{\X{}{} = (\X{}{1},\cdots,\X{}{I_C})\in \bigotimes_{i=1}^{I_C}\femspacecurvevector{i} \ \mid\ 
\X{}{\curveindex{k}{1}}(\curveindexrho{k}{1}) =
\X{}{\curveindex{k}{2}}(\curveindexrho{k}{2}) = 
\X{}{\curveindex{k}{3}}(\curveindexrho{k}{3})\ 1 \leq k \leq I_T \right\},
\]
where $\curveindexrho{k}{i} \in \{0,1\}$ encodes at which if its two endpoints
the discrete curve $\Gamma^m_{\curveindex{k}{i}}$ meets the 
$k$-th triple junction.

Given $\X{m}{}\in \underline{V}(\curve{m}{})$,
we find $(\chemicalDiscVec{m+1}{}, \X{m+1}{}, \kappa^{m+1}_\sigma)\in {\femspacebulkvectorsum}\times \underline{V}(\curve{m}{})\times V(\curve{m}{})$ such that the following conditions hold:\newline\newline
\noindent
\textbf{[Motion law]} For all ${\bv \varphi}\in{\femspacebulkvectorsum}$,
\begin{subequations} \label{eq:FEAmulti}
\begin{equation}\label{eq:FEAmulti-ML}
   \innerproduct{\nabla \chemicalDiscVec{m+1}{}}{\nabla\bv \varphi}{}{\Omega} 
   + \sum_{i=1}^{I_C}\sum_{j=1}^{I_P}\dcmap_{ji}\innerproduct{\pi^{h}_{i}\left[\frac{\X{m+1}{i} - \X{m}{i}}{\tau_m}\cdot\weightnormalnoindex{m}{i}\right]}{\varphi_{j}}{(h)}{\curve{m}{i}} = 0.
\end{equation}
\textbf{[Gibbs--Thomson law]} For all $\xi\in {V(\curve{m}{})}$,
\begin{equation}\label{eq:FEAmulti-GTL}
 \innerproduct{\kappa^{m+1}_\sigma}{\xi}{h}{\curve{m}{}}
       + \sum_{i=1}^{I_C}\sum_{j=1}^{I_P}\dcmap_{ji}\innerproduct{\chemicalDisc{m+1}{j}}{\xi_i}{(h)}{\curve{m}{i}}  = 0.
\end{equation}
\textbf{[Curvature vector]} For all $\vec{\eta}\in\underline{V}(\curve{m}{})$,
\begin{equation}\label{eq:FEAmulti-CV}
    \innerproduct{\kappa^{m+1}_\sigma\weightnormalnoindex{m}{}}{\vec{\eta}}{h}{\curve{m}{}} + \innerproduct{\sigma\sgrad\X{m+1}{}}{\sgrad\vec{\eta}}{}{\curve{m}{}} = 0.
\end{equation}
\end{subequations}

\begin{rem}[Linear system]
The linear system of equations arising at each time level of 
\eqref{eq:FEAmulti} is given by the obvious generalization of 
\eqref{eq:MatrixForm}, where the block matrix entries of \eqref{eq:MatrixForm}
are now defined by $A_\Omega = \operatorname{diag}{(A)_{j=1,\ldots,I_P}}$,
$C_\Gamma = \operatorname{diag}{(C_i)_{i=1,\ldots,I_C}}$
$\vec{N}_{\Omega,\Gamma} =
(\dcmap_{ji}\vec{N}_i)_{i=1,\ldots,I_C,j=1,\ldots,I_P}$,
$B_{\Omega,\Gamma} = (\dcmap_{ji}B_i)_{i=1,\ldots,I_C,j=1,\ldots,I_P}$
$\vec{D}_\Gamma = \operatorname{diag}{(\Vec{D}_i)_{i=1,\ldots,I_C}}$ and
$\metric{{E}_\Gamma} = \operatorname{diag}
{(\sigma_i\metric{E_i})_{i=1,\ldots,I_C}}$.
Once again the generalized system \eqref{eq:MatrixForm} can be reduced by
eliminating the final component $W^{m+1}_{I_P}$ from $\chemicalDiscVec{m+1}{}$.
We obtain the same block structure as in \eqref{eq:MatrixFormReduced},
with the new entries now given by
$\widehat A_\Omega = \operatorname{diag}{(A)_{j=1,\ldots,I_P-1}}$,
$\widehat {B}_{\Omega,\Gamma} = 
((\dcmap_{I_P,i}-\dcmap_{ji})B_i)_{i=1,\ldots,I_C,j=1,\ldots,I_P-1}$ and
$\widehat {N}_{\Omega,\Gamma} =
(\dcmap_{ji}\vec{N}_i)_{i=1,\ldots,I_C,j=1,\ldots,I_P-1}$.
\end{rem}

Finally, on using the techniques from Section~\ref{sec:ImprovementOfTheScheme},
we can adapt the approximation \eqref{eq:FEAmulti} to obtain a structure preserving
scheme that is unconditionally stable and that conserves the areas of the
enclosed phases exactly.
Given $(\X{m}{},\kappa^m_\sigma)\in \underline{V}(\curve{m}{})\times V(\curve{m}{})$, 
find $(\chemicalDiscVec{m+1}{}, \X{m+1}{},$\linebreak$ \kappa^{m+1}_\sigma)\in \femspacebulkvectorsum\times \underline{V}(\curve{m}{})\times V(\curve{m}{})$ such that
\newline\newline\noindent
\textbf{[Motion law]} For all $\bv \varphi\in \femspacebulkvectorsum$,
\begin{subequations} \label{eq:FEAnlmulti}
    \begin{equation}\label{eq:FEAnlmulti-ML}
        \innerproduct{\nabla \chemicalDiscVec{m+1}{}}{\nabla\bv \varphi}{}{\Omega} 
        + \sum_{i=1}^{I_C}\sum_{j=1}^{I_P}\dcmap_{ji}\innerproduct{\pi^h_i\left[\frac{\X{m+1}{i} - \X{m}{i}}{\tau_m}\cdot\weightnormalnoindex{m+\frac{1}{2}}{i}\right]}{\varphi_{j}}{(h)}{\curve{m}{i}} = 0.
    \end{equation}
    \textbf{[Gibbs--Thomson law]} For all $\xi\in V(\curve{m}{})$,
    \begin{equation}\label{eq:FEAnlmulti-GTL}
 \innerproduct{\kappa^{m+1}_\sigma}{\xi}{h}{\curve{m}{}}
       + \sum_{i=1}^{I_C}\sum_{j=1}^{I_P}\dcmap_{ji}\innerproduct{\chemicalDisc{m+1}{j}}{\xi_i}{(h)}{\curve{m}{i}}  = 0.
    \end{equation}
    \textbf{[Curvature vector]} For all $\vec{\eta}\in\underline{V}(\curve{m}{})$,
    \begin{equation}\label{eq:FEAnlmulti-CV}
        \innerproduct{\kappa^{m+1}_\sigma\weightnormalnoindex{m+\frac{1}{2}}{}}{\vec{\eta}}{h}{\curve{m}{}} + \innerproduct{\sigma\sgrad\X{m+1}{}}{\sgrad\vec{\eta}}{}{\curve{m}{}} = 0.
    \end{equation}
\end{subequations}

We conclude this section by stating theoretical results for the generalized
schemes. Their proofs are straightforward adaptations of the proofs
of Theorems~\ref{thm:exist_unique},
\ref{thm:StabilityLinear} and \ref{thm:AP_nonlinear}.

\begin{asm} \label{asm:4}
Assume that 
$\operatorname{span}{\{\weightnormal{m}{i}{j}\}_{1\leq j\leq N_i-1}} \not= 
\{\vec0\}$ for $1 \leq i \leq I_C$ and
\[
        \operatorname{span}\left\{
\sum_{i=1}^{I_C}\sum_{j=1}^{I_P}\dcmap_{ji}\innerproduct{\weightnormalnoindex{m}{i}}{\varphi_{j}}{(h)}{\curve{m}{i}}
 \mid\ \bv \varphi\in\femspacebulkvectorsum\right\} = \bR^2.
\]
\end{asm}

\begin{thm}
    Suppose that Assumption~\ref{asm:4} holds, and that $m\geq0$.
    Then, there exists a unique solution
$(\chemicalDiscVec{m+1}{},\X{m+1}{},\kappa^{m+1}_\sigma)\in\femspacebulkvectorsum\times\femspacegammavector\times\femspacegamma$ to \eqref{eq:FEAmulti}.
Moreover, any solution to \eqref{eq:FEAmulti} or \eqref{eq:FEAnlmulti}
satisfies the stability bound
    \begin{equation*}
        |\curve{m+1}{}|_\sigma + \tau_m\|\nabla\chemicalDiscVec{m+1}{}\|^2_2\leq |\curve{m}{}|_\sigma.
    \end{equation*}
    Finally, a solution to \eqref{eq:FEAnlmulti} satisfies
    \begin{equation*}
        |\Omega^{m+1}_j| = |\Omega^m_j| \quad 1\leq j \leq I_P.
    \end{equation*}
\end{thm}

\section{Numerical results}\label{sec:NumericalResults}
\newcommand{\errorXx}{\|\Gamma^h - \Gamma\|_{L^\infty}}
\newcommand{\errorUu}{\|\bv W^h - \bv w\|_{L^\infty}}
\newcommand{\dist}{\operatorname{dist}}

We implemented the fully discrete finite element approximations 
\eqref{eq:FEAmulti} and \eqref{eq:FEAnlmulti} within the
finite element toolbox ALBERTA, see \cite{Alberta}. The arising 
linear systems of the form \eqref{eq:MatrixFormReduced} 
are solved with a GMRes iterative solver, applying as preconditioner
a least squares solution of the block matrix in 
\eqref{eq:MatrixFormReduced} without the projection matrices
$\metric{P}$. For the computation of the least squares solution we employ
the sparse factorization package SPQR, see \cite{Davis11}.

For the triangulation $\mathfrak{T}^m$ of the bulk domain $\Omega$, that is used
for the bulk finite element space $S^m$, we use an adaptive mesh that uses 
fine elements close to the interface $\Gamma^m$ and coarser elements away 
from it. The precise strategy is as described in 
\cite[\S5.1]{BarrettGarckeRobert2010} 
and for a domain $\Omega=(-H,H)^d$ and two integer parameters
$N_c < N_f$ results in elements with maximal diameter approximately equal to
$h_f = \frac{2H}{N_f}$ close to $\Gamma^m$ and elements with maximal diameter
approximately equal to $h_c = \frac{2H}{N_c}$ far away from it. 
For all our computations we use $H=4$.
An example adaptive mesh is shown in Figure~\ref{fig:2d_db_plus_one}, below.

We stress that due to the unfitted nature of our finite element approximations,
special quadrature rules need to be employed in order to assemble terms that
feature both bulk and surface finite element functions. For all the
computations presented in this section, we use true integration for these
terms, and we refer to
\cite{BarrettGarckeRobert2010,Nurnberg202203} for details on the practical
implementation.
Throughout this section we use (almost) uniform time steps, in that
$\tau_m=\tau$ for $m=0,\ldots, M-2$ and 
$\tau_{M-1} = T - t_{m-1} \leq \tau$. 
Unless otherwise stated, we set $\sigma=1$.

\subsection{Convergence experiment}\label{sec:ConvergenceExperiments}
\begin{figure}[H]
        \centering
        \includegraphics[keepaspectratio, scale=0.2]{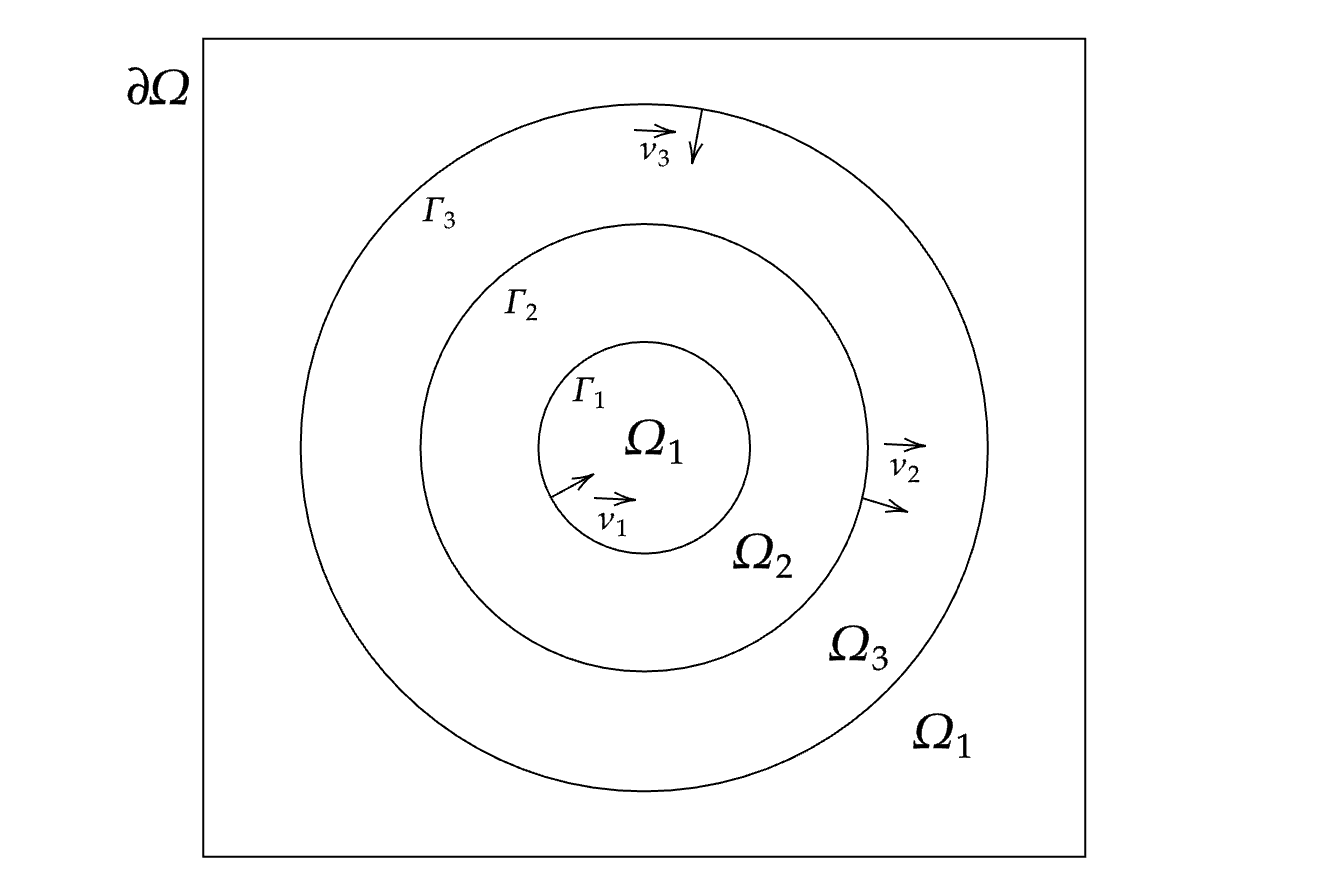}
        \caption{A network of curves which consists of three concentric circles.}\label{fig:3concentriccircles}
    \end{figure}

    In order to validate our proposed schemes, we utilize the following exact
solution for a network of three concentric circles.
    Let $0<R_1(t)<R_2(t)<R_3(t)$ and $\curve{}{i}(t) := \partial B(0,R_i(t))$ 
for each $i\in\naturalset{3}$.
    Suppose that $\normalnoindex{}{1}$, $-\normalnoindex{}{2}$ and 
$\normalnoindex{}{3}$ are pointing towards the origin, and 
let $\Omega_1(t) = B(0,R_1(t)) \cup (\Omega \setminus \overline{B(0,R_3(t))})$, 
$\Omega_j(t) = B(0,R_j(t)) \setminus \closure{B(0,R_{j-1}(t))}$, $j=2,3$.
Hence, with the notation from \S\ref{sec:ProblemSettingMulti} we have
$I_C = 3$, $I_P = 3$, $I_T = 0$
and
$\dcmap = \begin{pmatrix}-1 &0&1 \\ 1 & 1 & 0 \\ 0& -1 &-1 \end{pmatrix}$.
    See Figure \ref{fig:3concentriccircles} for the setting. 

We shall prove in Appendix~\ref{sec:B} that 
$(\bv w,\{\curve{}{}(t)\}_{0\leq t\leq T})$ is a solution to
\eqref{eq:StrongForm_Multi} with $\sigma=1$ 
if the three radii satisfy the differential
algebraic equations{:}
\begin{subequations} \label{eq:ODE}
    \begin{equation}\label{eq:3circleRadiusFormula}
            R'_2(t) = -F(R_2(t)),\quad
            R_1(t) = \sqrt{R_2(t)^2 - A_2},\quad
            R_3(t) = \sqrt{R_2(t)^2 + A_3},
    \end{equation}
    where 
    $F(u) := (\frac{1}{\sqrt{u^2 - A_2}} + \frac{1}{u} + \frac{1}{\sqrt{u^2 + A_3}})/({u\log{\frac{u^2 + A_3}{u^2 - A_2}}})$ for $u\in(\sqrt{A_2},\infty)$,
$A_2 := R_2(0)^2 - R_1(0)^2$, and $A_3 := R_3(0)^2 - R_2(0)^2$,
and if the three chemical potentials are given by
    \begin{align}\label{eq:3circleChemicalFormula}
                    w_1(x,\cdot) &= -\frac{1}{3R_2}-\frac{2}{3R_3} + \begin{cases}
                        \frac{1}{3}\alpha\log{\frac{R_3^4}{R_1^3R_2}}, & |x|\leq R_1\\
                        - \frac{1}{3}\alpha\log{\frac{R_2|x|^3}{R_3^4}}, & R_1<|x|\leq R_3\\
                        \frac{1}{3}\alpha\log{\frac{R_3}{R_2}}, & R_3 < |x|,
                    \end{cases} \nonumber \\
                    w_2(x,\cdot) &= \begin{cases}
                        -\frac{1}{R_1}-\frac{1}{3R_2}-\frac{2}{3R_3}+ \frac{1}{3}\alpha\log{\frac{R_3^4}{R_1^3R_2}}, & |x|\leq R_1\\
                        -\frac{1}{R_1}-\frac{1}{3R_2}-\frac{2}{3R_3} + \frac{1}{3}\alpha\log{\frac{R_3^4|x|^3}{R_1^6R_2}}, & R_1 < |x| \leq R_2\\
                        \frac{2}{3R_2}+\frac{1}{3R_3}+\frac{1}{3}\alpha\log{\frac{R_2^2}{R_3^2}}, & R_2 < |x|,
                    \end{cases} \nonumber \\
                    w_3(x,\cdot) &= \begin{cases}
                        -\frac{1}{R_1}-\frac{4}{3R_2}-\frac{2}{3R_3}+\frac{1}{3}\alpha\log{\frac{R_2^2R_3^4}{R_1^6}}, & |x| < R_2\\
                        -\frac{1}{3R_2}+\frac{1}{3R_3}+\frac{1}{3}\alpha\log{\frac{|x|^3}{R_2R_3^2}}, & R_2\leq |x|<R_3\\
                        -\frac{1}{3R_2}+\frac{1}{3R_3}+\frac{1}{3}\alpha\log{\frac{R_3}{R_2}}, & R_3\leq |x|,
                    \end{cases}
    \end{align}
\end{subequations}
    where $\alpha(t) := \constalpha{(t)}$.

In order to accurately compute the radius $R_2(t)$, rather than numerically
solving the ODE in \eqref{eq:3circleRadiusFormula}, 
we employ a root finding algorithm for the equation
        \begin{equation}\label{eq:RootFindingEquation}
            0 = t + \int_{R_2(0)}^{R_2(t)}\frac{1}{F(u)}\,du,
        \end{equation}
following similar ideas in \cite{BarrettGarckeRobert2010,Nurnberg202203}.

For the initial radii $R_1(0) = 2$, $R_2(0) = 2.5$, $R_3(0) = 3$ and 
the time interval $[0,T]$ with $T=\frac12$, 
so that $R_1(T) \approx 1.60$, $R_2(t) \approx 2.20$ and 
$R_3(T) \approx 2.75$, we perform a convergence 
experiment for the true solution \eqref{eq:ODE}.
To this end, for $i=0\to 4$, we set 
$N_f = \frac12 K = 2^{7+i}$, $N_c = 4^{i}$
and $\tau= 4^{3-i}\times10^{-3}$. 
In Tables~\ref{tab:DMS2dT05} and \ref{tab:fdMS2dT05}
we display the errors 
\[\errorXx = \max_{m=1,\ldots, M} 
\max_{i=1,\ldots,3} \max_{j=1,\ldots,N_i} \dist(\vec{q}^m_{i,j}, \Gamma_i(t_m))
\]
and
\[
\errorUu = \max_{m=1,\ldots, M} \max_{i=1,\ldots,3} \|W^m_i - I^m w_i(\cdot,t_m)\|_{L^\infty(\Omega)},
\]
where $I^m : C^0(\overline\Omega) \to S^m$ denotes the standard 
interpolation operator,
for the schemes \eqref{eq:FEAmulti} and \eqref{eq:FEAnlmulti},
respectively.
We also let $K_\Omega^m$ denote the number of
degrees of freedom of $S^m$, and define $h^m_\Gamma = \max_{i=1,\ldots,3}
\max_{j = 1,\ldots,N_i} \sigma^m_{i,j}$,
as well as $v_\Delta^M = \max_{j=1,\ldots,I_P} |\,|\Omega^M_j| -
|\Omega^0_j| \,|$.
As expected, we observe true volume preservation for the scheme 
\eqref{eq:FEAnlmulti} in Table~\ref{tab:fdMS2dT05}, 
up to solver tolerance, while the relative volume
loss in Table~\ref{tab:DMS2dT05} decreases as $\tau$ becomes smaller. 
Surprisingly, the two error quantities $\errorXx$ and $\errorUu$ are generally
lower in Table~\ref{tab:DMS2dT05} compared to Table~\ref{tab:fdMS2dT05},
although the difference becomes smaller with smaller discretization parameters.

\begin{table}[h]
\center
\begin{tabular}{c|c|c|c|c|c|c|c}
 $h_{f}$ & $h^M_\Gamma$ & $\tau$ & $\errorUu$ & $\errorXx$ & $K^M_\Omega$ & $N$
 & $v_\Delta^M$ \\ \hline
6.2500e-02 &1.3487e-01 &6.4e-2&3.4879e-02 &4.7230e-03 & 3633 & 384 & 3.8e-02\\
3.1250e-02 &6.7465e-02 &1.6e-2&1.2459e-02 &4.3450e-03 & 7321 & 768 & 9.7e-03\\
1.5625e-02 &3.3747e-02 &4.0e-3&4.8323e-03 &2.8582e-03 & 14881 & 1536& 2.4e-03\\
7.8125e-03 &1.6878e-02 &1.0e-3&2.0917e-03 &1.6382e-03 & 30161 & 3072 & 6.1e-04\\
3.9062e-03 &8.4405e-03 &2.5e-4&1.4139e-03 &8.4224e-04 & 72537 & 6144 & 1.5e-04\\
\end{tabular}
\caption{Convergence test for \eqref{eq:ODE} over the time interval 
$[0,\frac12]$ for the scheme \eqref{eq:FEAmulti}.}
\label{tab:DMS2dT05}
\end{table}%
\begin{table}[h]
\center
\begin{tabular}{c|c|c|c|c|c|c|c}
 $h_{f}$ & $h^M_\Gamma$ & $\tau$ & $\errorUu$ & $\errorXx$ & $K^M_\Omega$ & $N$
 & $v_\Delta^M$ \\ \hline
6.2500e-02 &1.3469e-01&6.4e-2 &3.8031e-02 &1.3628e-02 &3605 &384 &$<10^{-10}$\\
3.1250e-02 &6.7442e-02&1.6e-2 &1.3442e-02 &6.7806e-03 &7285 &768 &$<10^{-10}$\\
1.5625e-02 &3.3744e-02&4.0e-3 &5.0754e-03 &3.4721e-03 &14905&1536&$<10^{-10}$\\
7.8125e-03 &1.6878e-02&1.0e-3 &2.1316e-03 &1.7904e-03 &30193&3072&$<10^{-10}$\\
3.9062e-03 &8.4404e-03&2.5e-4 &1.4242e-03 &8.8016e-04 &72537&6144&$<10^{-10}$\\
\end{tabular}
\caption{Convergence test for \eqref{eq:ODE} over the time interval 
$[0,\frac12]$ for the scheme \eqref{eq:FEAnlmulti}.}
\label{tab:fdMS2dT05}
\end{table}%

For all the following numerical simulations, we always employ the fully
structure-preserving scheme \eqref{eq:FEAnlmulti}.

\subsection{Evolutions with equal surface energies} \label{sec:83}
\subsubsection{3 phases}
In the next set of experiments, we investigate how a standard double bubble and
a disk evolve, when one phase is made up of the left bubble, and the other
phase is made up of the right bubble and the disk.
With the notation from \S\ref{sec:ProblemSettingMulti} we have
$I_C = 4$, $I_P = 3$, $I_T = 2$, 
$(\curveindex{1}{1},\curveindex{1}{2},\curveindex{1}{3}) = (\curveindex{2}{1},\curveindex{2}{2},\curveindex{2}{3}) = (1,2,3)$
and
$\dcmap = \begin{pmatrix}0 &-1&1&0 \\ 1 & 0 &-1&-1 \\ -1 &1&0&1 \end{pmatrix}$.

\begin{figure}
\center
\includegraphics[angle=-90,width=0.18\textwidth]{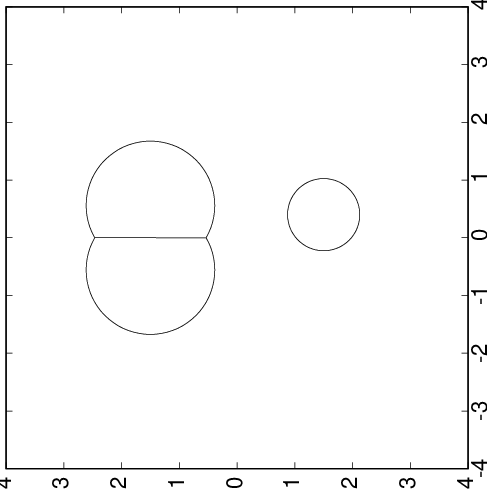}
\includegraphics[angle=-90,width=0.18\textwidth]{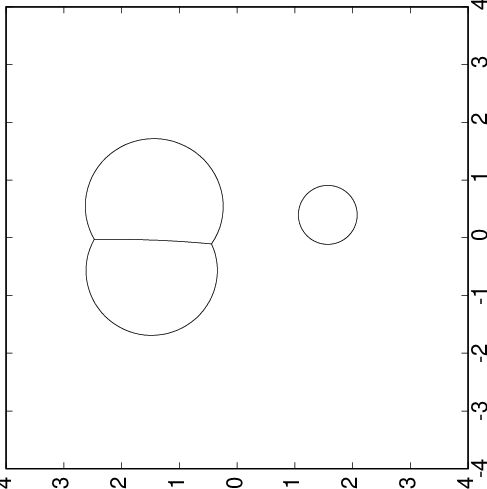}
\includegraphics[angle=-90,width=0.18\textwidth]{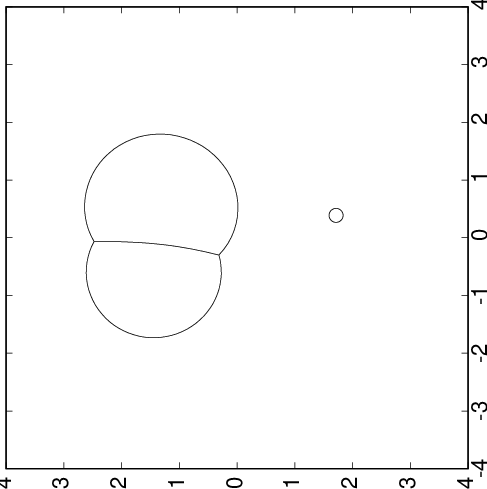}
\includegraphics[angle=-90,width=0.18\textwidth]{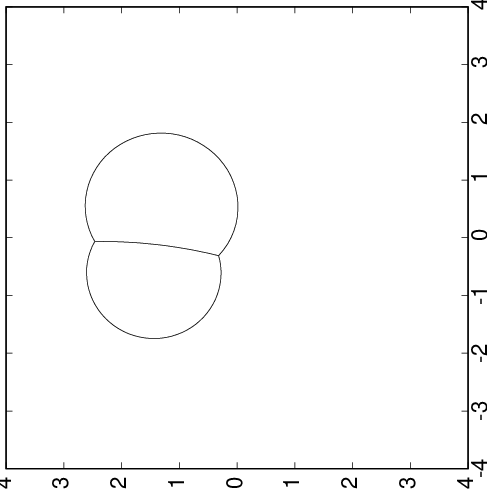}
\includegraphics[angle=-90,width=0.18\textwidth]{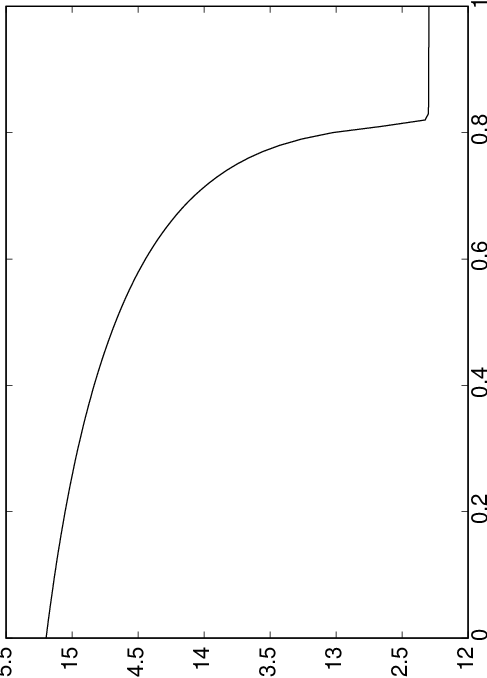} 
\\[5mm]
\includegraphics[angle=-0,width=0.18\textwidth]{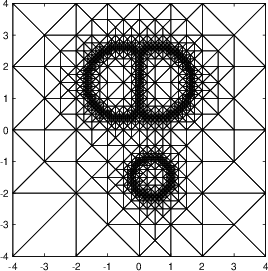}
\includegraphics[angle=-0,width=0.18\textwidth]{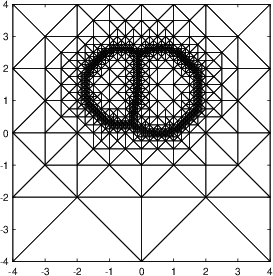}
\caption{The solution at times $t=0, 0.4, 0.8, 1$, and a plot of the discrete 
energy over time. Below we show the adaptive bulk mesh at times $t=0$ and
$t=1$.
}
\label{fig:2d_db_plus_one}
\end{figure}

\begin{figure}
\center
\includegraphics[angle=-90,width=0.18\textwidth]{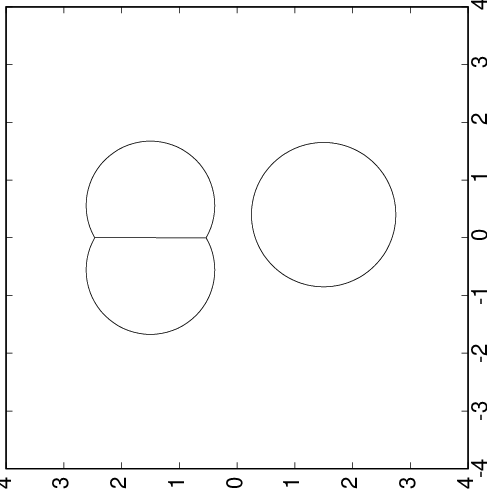}
\includegraphics[angle=-90,width=0.18\textwidth]{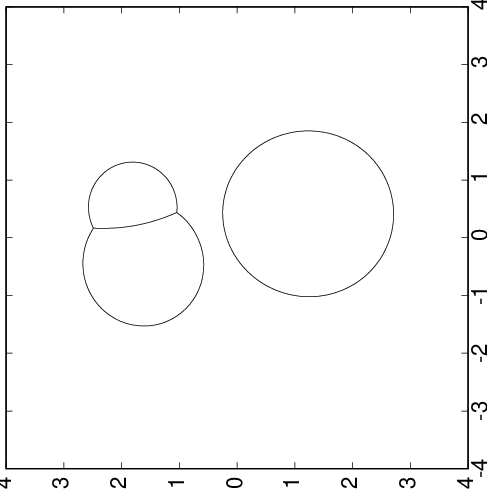}
\includegraphics[angle=-90,width=0.18\textwidth]{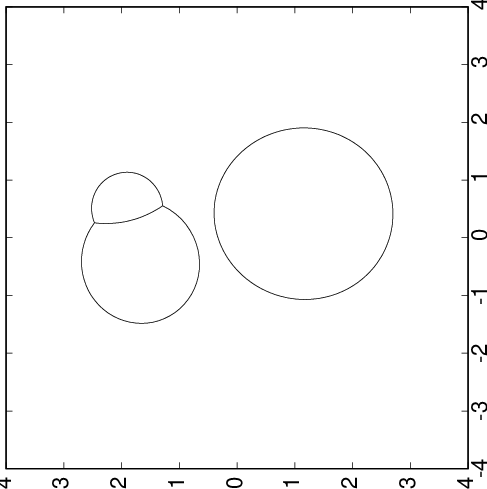}
\includegraphics[angle=-90,width=0.18\textwidth]{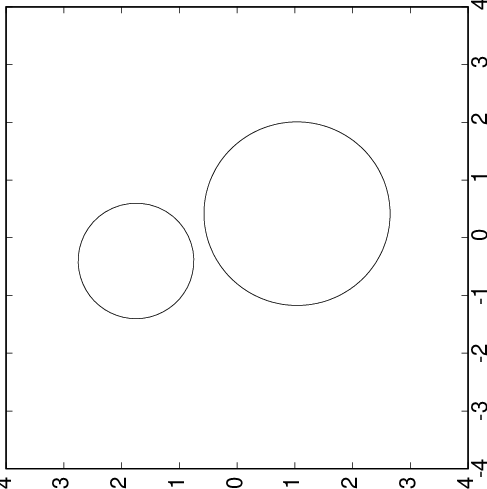}
\includegraphics[angle=-90,width=0.18\textwidth]{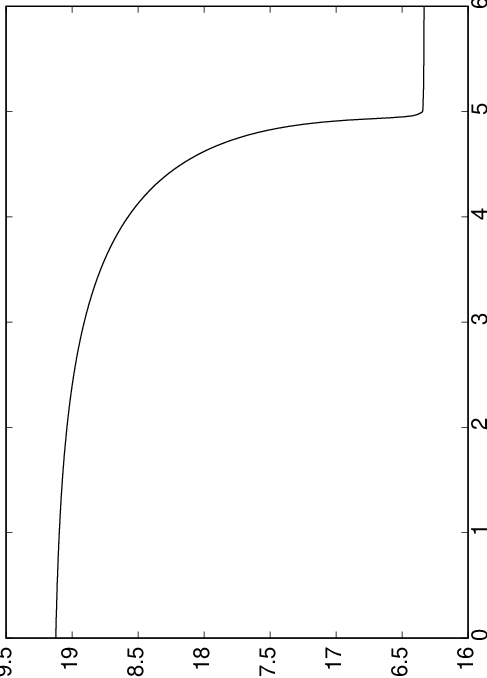}
\caption{The solution at times $t=0, 4, 4.5, 6$, and a plot of the discrete 
energy over time.
}
\label{fig:2d_db_plus_bigone}
\end{figure}

\begin{figure}
\center
\includegraphics[angle=-90,width=0.18\textwidth]{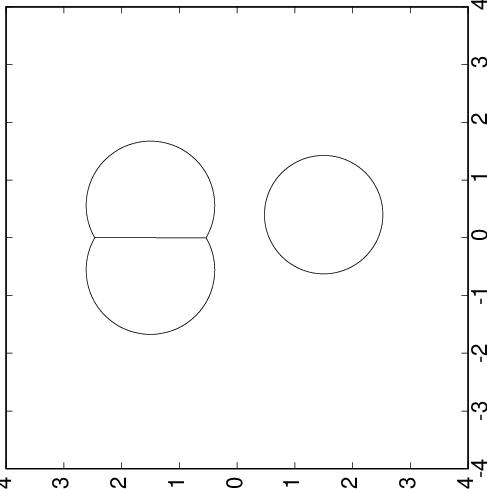}
\includegraphics[angle=-90,width=0.18\textwidth]{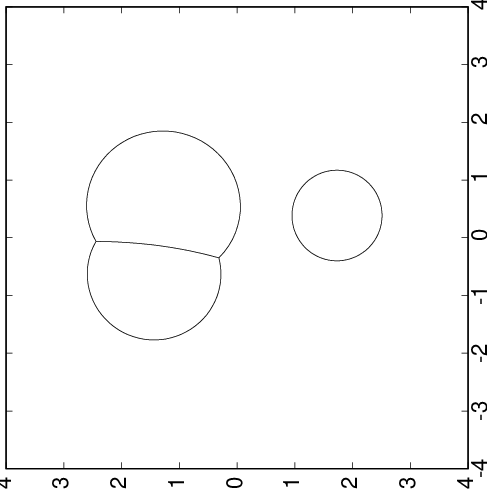}
\includegraphics[angle=-90,width=0.18\textwidth]{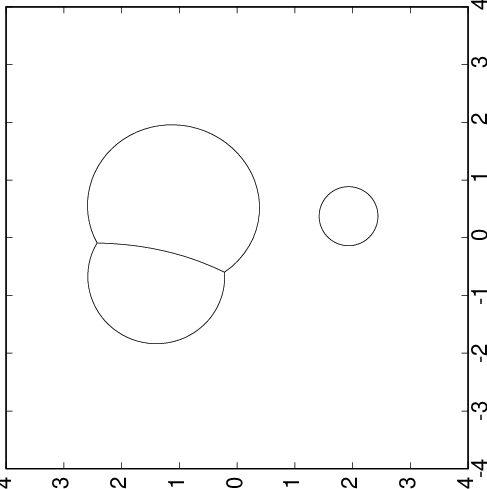}
\includegraphics[angle=-90,width=0.18\textwidth]{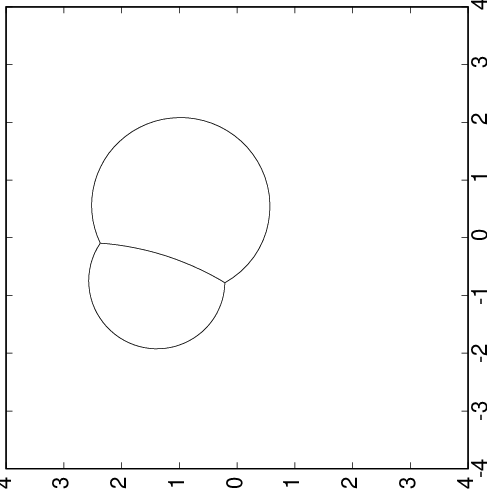}
\includegraphics[angle=-90,width=0.18\textwidth]{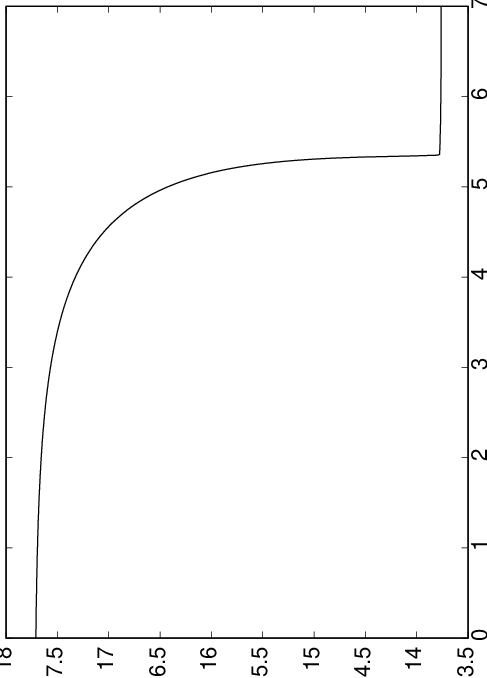}
\caption{The solution at times $t=0, 4, 5, 7$, and a plot of the discrete 
energy over time.
}
\label{fig:2d_db_plus_04one}
\end{figure}

The two bubbles of the double bubble enclose an area of about $3.139$ each,
while the disk has an initial radius of $\frac58$, meaning it initially
encloses an area of $\frac{25\pi}{64} \approx 1.227$. During the evolution the
disk vanishes, and the right bubble grows correspondingly, see 
Figure~\ref{fig:2d_db_plus_one}. We note that our theoretically framework
does not allow for changes of topology, e.g., the vanishing of curves. Hence, in
our computations we perform heuristic surgeries whenever a curve becomes too
short.
Here a closed curve is simply discarded, while a curve that was part of a 
network is removed. This will leave two triple junctions, where only two 
curves meet, and the involved curves can be glued together so that the
simulation can continue.
Repeating the simulation with a bigger initial disk gives the results in
Figure~\ref{fig:2d_db_plus_bigone}. Here the radius is $\frac54$, so that the
enclosed area is $4.909$. Now the disk grows at the expense of the right
bubble, so that eventually two separate phases remain.
With the next simulation we demonstrate that in the given setup, which of
the two components of phase 1 survives is not down the initial size. 
In particular, we
allow the initial disk to have area $3.300$, so that it is bigger than the
other component of the same phase: the right bubble in the double bubble. And
yet, due to the perimeter of the bubble being overall cheaper than the boundary
of the disk, the latter shrinks to extinction. See 
Figure~\ref{fig:2d_db_plus_04one}.

\begin{figure}
\center
\includegraphics[angle=-90,width=0.18\textwidth]{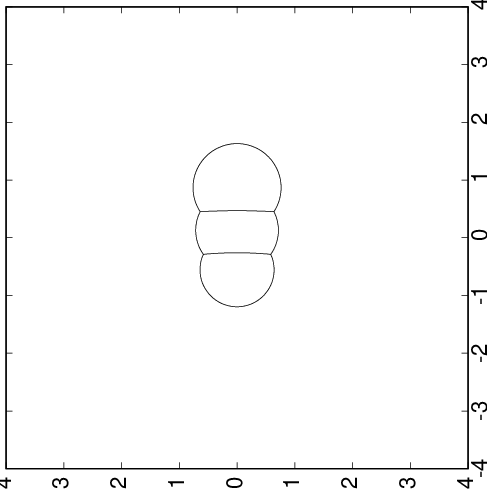}
\includegraphics[angle=-90,width=0.18\textwidth]{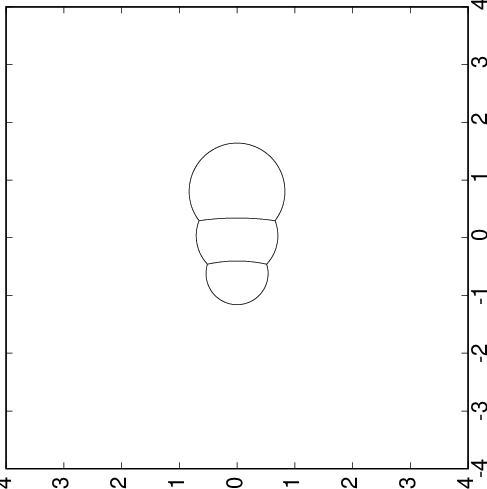}
\includegraphics[angle=-90,width=0.18\textwidth]{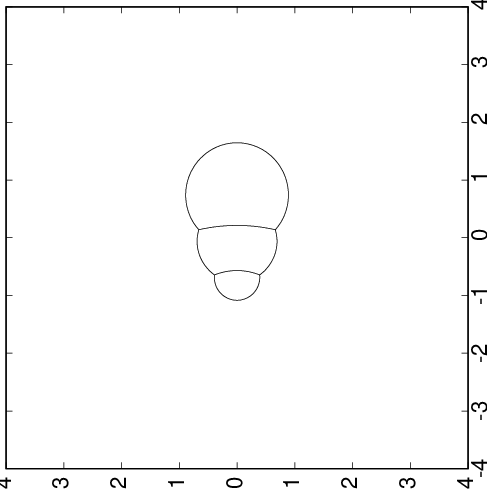}
\includegraphics[angle=-90,width=0.18\textwidth]{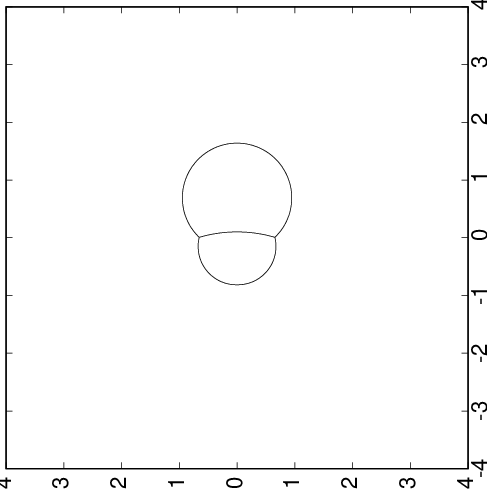}
\includegraphics[angle=-90,width=0.18\textwidth]{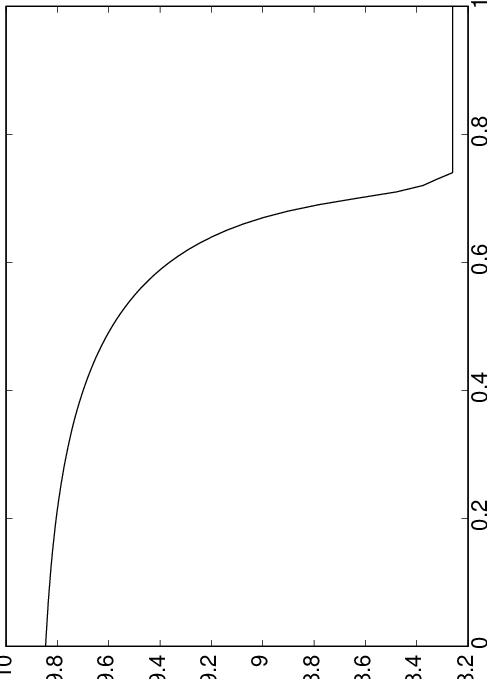}
\caption{The solution at times $t=0, 0.4, 0.6, 1$, and a plot of the discrete 
energy over time.
}
\label{fig:2d_tb_NoP3}
\end{figure}

For the next experiment we start from a nonstandard triple bubble, where we
choose the right most bubble to have area $1.5$, while the other two bubbles
have unit area. We assign the two outer bubbles to belong to the same phase.
In particular, with the notation from \S\ref{sec:ProblemSettingMulti} we have
$I_C = 6$, $I_P = 3$, $I_T = 4$, 
$(\curveindex{1}{1},\curveindex{1}{2},\curveindex{1}{3}) = (1,2,5)$,
$(\curveindex{2}{1},\curveindex{2}{2},\curveindex{2}{3}) = (1,3,5)$,
$(\curveindex{3}{1},\curveindex{3}{2},\curveindex{3}{3}) = (2,4,6)$,
$(\curveindex{4}{1},\curveindex{4}{2},\curveindex{4}{3}) = (3,4,6)$
and
\begin{equation*}
\dcmap = \begin{pmatrix} -1 & 0 & 0 & -1 & 1 & 1\\ 0 & 1 & 1 & 0 & -1 & -1\\ 1 & -1 & -1 & 1 & 0 & 0 \end{pmatrix}.
\end{equation*}
We observe that during the evolution the larger bubble on the right grows at the expense
of the left bubble, until the latter one vanishes completely. The remaining
interfaces then evolve towards a standard double bubble with enclosed areas 1
and 2.5. See 
Figure~\ref{fig:2d_tb_NoP3}.

\begin{figure}
\center
\includegraphics[angle=-90,width=0.18\textwidth]{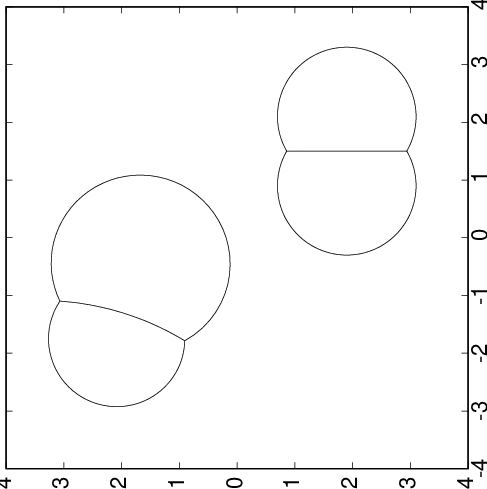}
\includegraphics[angle=-90,width=0.18\textwidth]{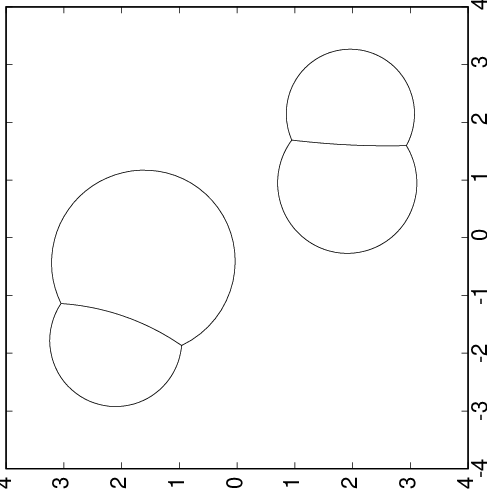}
\includegraphics[angle=-90,width=0.18\textwidth]{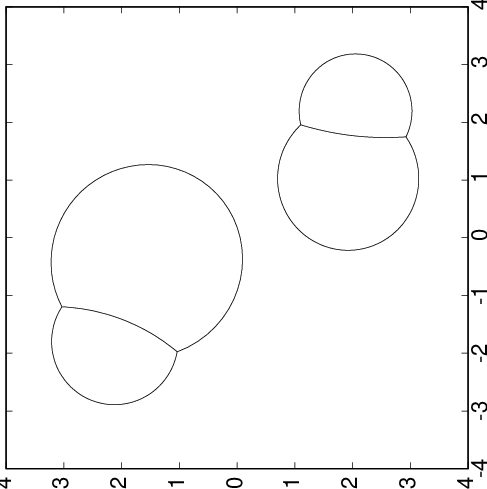}
\includegraphics[angle=-90,width=0.18\textwidth]{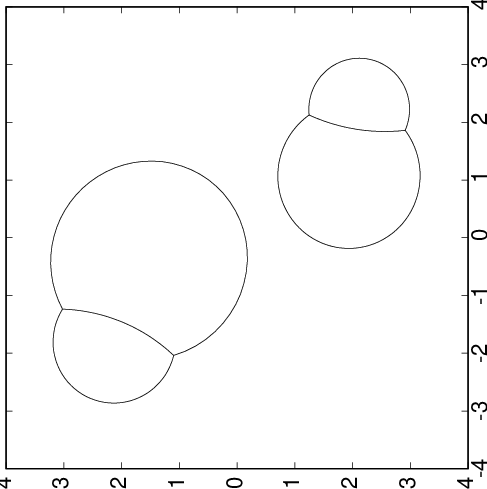}
\includegraphics[angle=-90,width=0.18\textwidth]{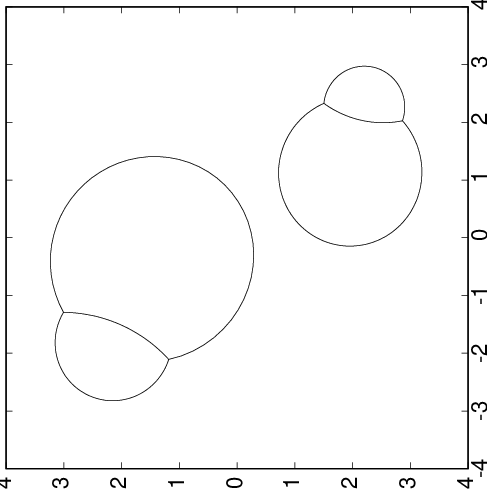}
\includegraphics[angle=-90,width=0.18\textwidth]{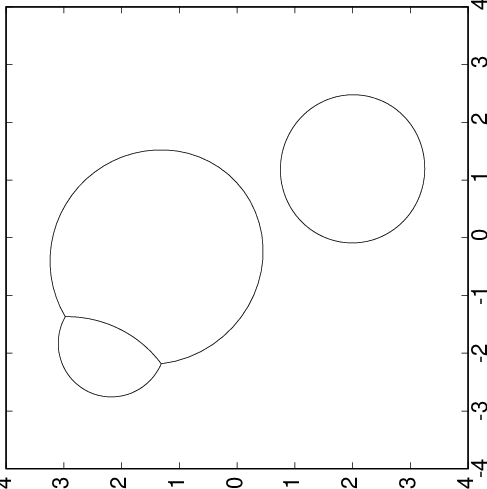}
\includegraphics[angle=-90,width=0.18\textwidth]{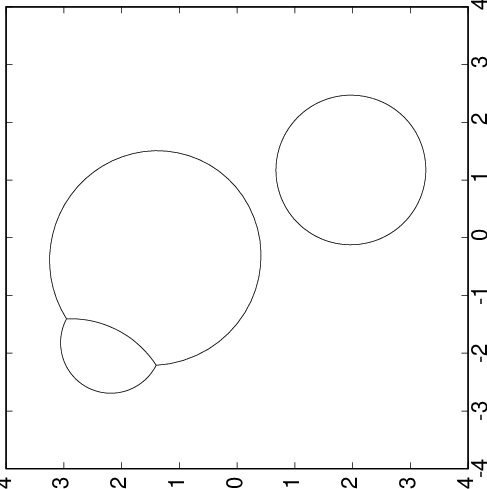}
\includegraphics[angle=-90,width=0.18\textwidth]{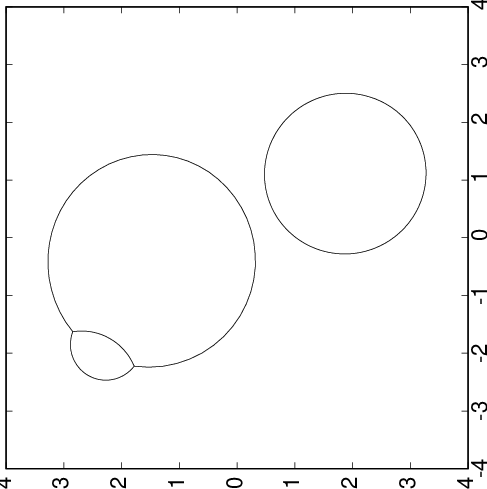}
\includegraphics[angle=-90,width=0.18\textwidth]{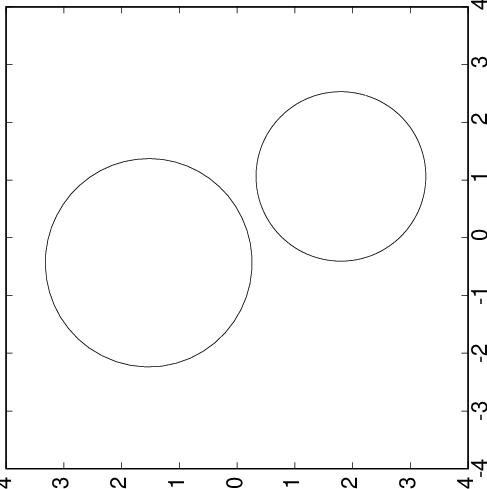}
\includegraphics[angle=-90,width=0.18\textwidth]{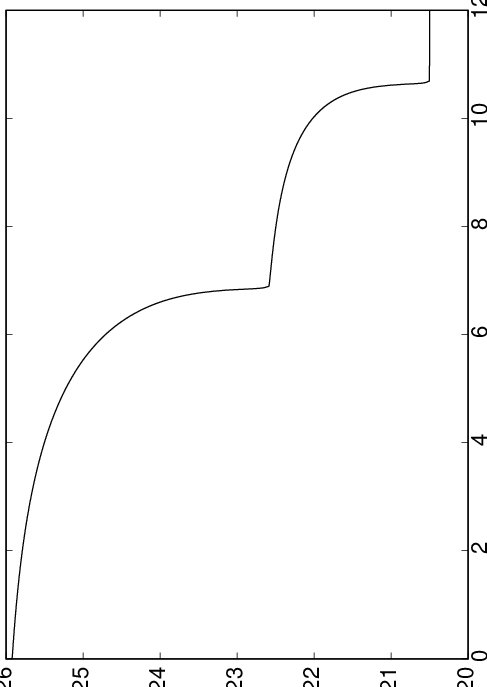}
\caption{The solution at times $t=0, 2, 4, 5, 6, 7, 8, 10, 12$, 
and a plot of the discrete energy over time.
}
\label{fig:2d_db_times2}
\end{figure}

In the final numerical simulation for the setting with 3 phases, we consider
the evolution of two double bubbles. 
With the notation from \S\ref{sec:ProblemSettingMulti} we have
$I_C = 6$, $I_P = 3$, $I_T = 4$, 
$(\curveindex{1}{1},\curveindex{1}{2},\curveindex{1}{3}) = 
(\curveindex{2}{1},\curveindex{2}{2},\curveindex{2}{3}) = (1,2,3)$,
$(\curveindex{3}{1},\curveindex{3}{2},\curveindex{3}{3}) = 
(\curveindex{4}{1},\curveindex{4}{2},\curveindex{4}{3}) = (4,5,6)$
and
\begin{equation*}
\dcmap = \begin{pmatrix} 
0 & -1 & 1 & 0 & -1 & 1\\
1 & 0 & -1 & 1 & 0 & -1\\
-1 & 1 & 0 & -1 & 1 & 0
\end{pmatrix}.
\end{equation*}
The first bubble is chosen with enclosing
areas $3.14$ and $6.48$, 
while the second double bubbles encloses two areas of size $3.64$. In each
case, the left bubble is assigned to phase 1, while the right bubbles are
assigned to phase 2. In this way, the lower double bubble holds the larger
portion of phase 1, while the upper double bubble holds the larger portion of
phase 2. Consequently, each double bubble evolves to a single disk that
contains just one phase. See
Figure~\ref{fig:2d_db_times2}.

\subsubsection{4 phases}

In the next set of experiments, we investigate simulations 
for a nonstandard triple bubble, with one of the bubbles making a phase
with a separate disk.
In particular, with the notation from \S\ref{sec:ProblemSettingMulti} we have
$I_C = 7$, $I_P = 4$, $I_T = 4$, 
$(\curveindex{1}{1},\curveindex{1}{2},\curveindex{1}{3}) = (1,2,5)$,
$(\curveindex{2}{1},\curveindex{2}{2},\curveindex{2}{3}) = (1,3,5)$,
$(\curveindex{3}{1},\curveindex{3}{2},\curveindex{3}{3}) = (2,4,6)$,
$(\curveindex{4}{1},\curveindex{4}{2},\curveindex{4}{3}) = (3,4,6)$
and
\begin{equation*}
\dcmap = \begin{pmatrix} 
    0 & 0 & 0 & -1 & 0 & 1 & -1\\
    0 & 1 & 1 & 0 & -1 & -1 & 0\\
    -1 & 0 & 0 & 0 & 1 & 0 & 0\\
    1 & -1 & -1 & 1 & 0 & 0 & 1
\end{pmatrix}.
\end{equation*}

\begin{figure}
\center
\includegraphics[angle=-90,width=0.18\textwidth]{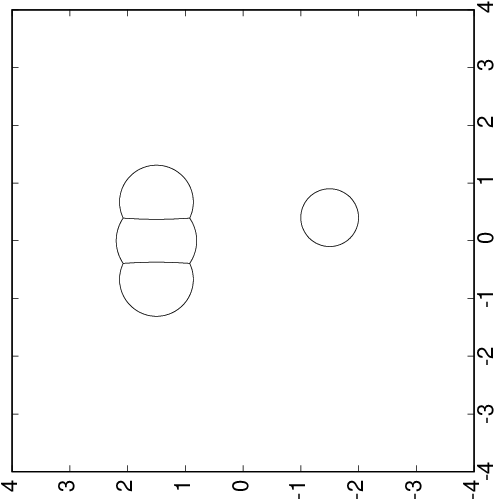}
\includegraphics[angle=-90,width=0.18\textwidth]{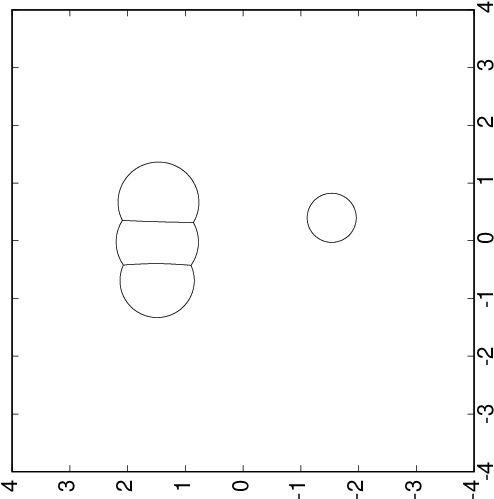}
\includegraphics[angle=-90,width=0.18\textwidth]{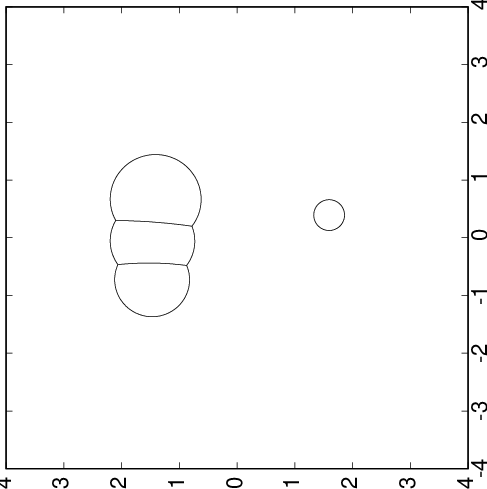}
\includegraphics[angle=-90,width=0.18\textwidth]{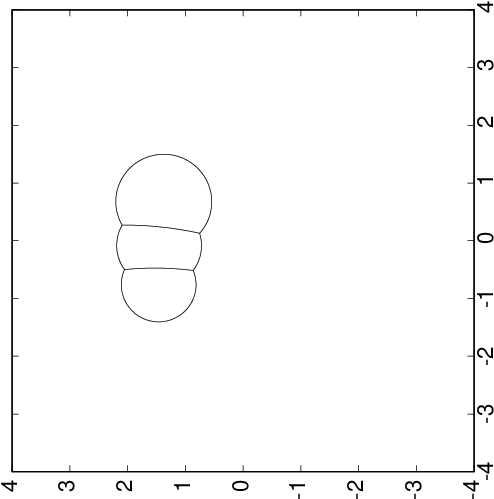}
\includegraphics[angle=-90,width=0.18\textwidth]{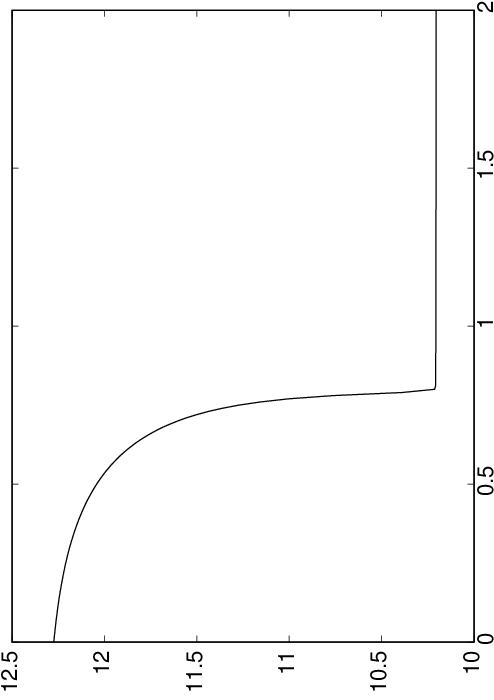}
\caption{The solution at times $t=0, 0.4, 0.7, 1$, and a plot of the discrete 
energy over time.
}
\label{fig:2d_tb_unstable_plus_one}
\end{figure}

\begin{figure}
\center
\includegraphics[angle=-90,width=0.18\textwidth]{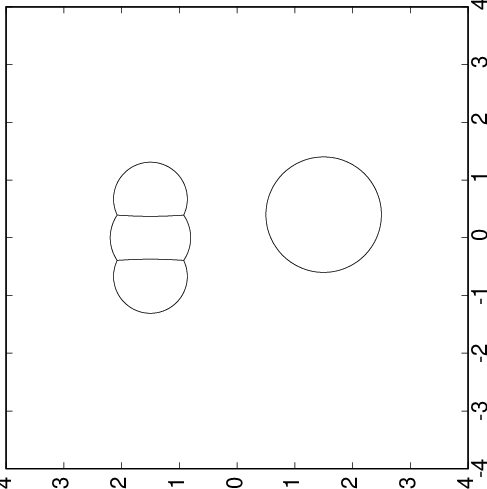}
\includegraphics[angle=-90,width=0.18\textwidth]{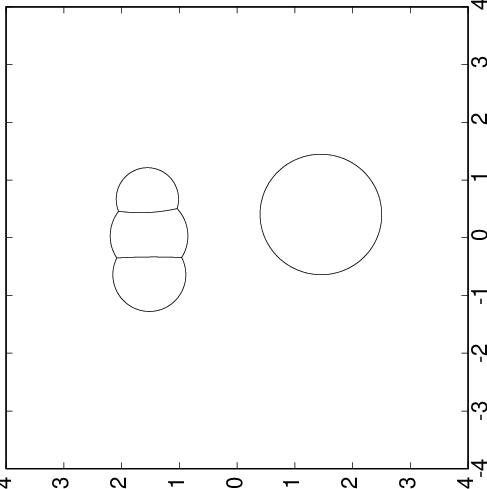}
\includegraphics[angle=-90,width=0.18\textwidth]{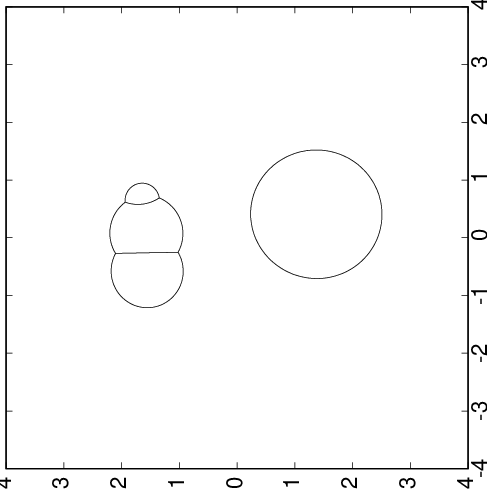}
\includegraphics[angle=-90,width=0.18\textwidth]{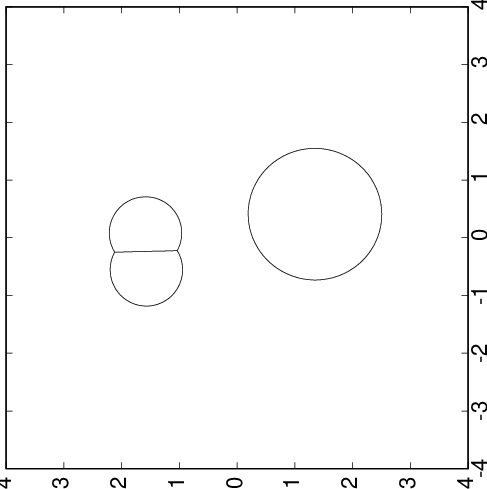}
\includegraphics[angle=-90,width=0.18\textwidth]{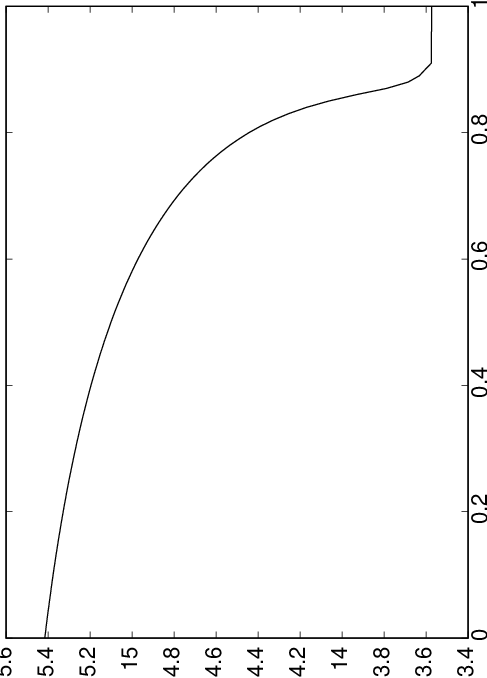}
\caption{The solution at times $t=0, 0.4, 0.8, 1$, and a plot of the discrete 
energy over time.
}
\label{fig:2d_tb_unstable_plus_bigone}
\end{figure}

The three bubbles of the triple bubble enclose an area of unity each,
while the disk has an initial radius of $\frac12$, meaning it initially
encloses an area of $\frac{\pi}{4} \approx 0.785$. During the evolution the
disk vanishes, and the right bubble grows correspondingly, see 
Figure~\ref{fig:2d_tb_unstable_plus_one}.
Repeating the experiment with initial data where the disk is a unit disk leads
to the evolution in Figure~\ref{fig:2d_tb_unstable_plus_bigone},
where the disk now expands and survives, at the expense of the right bubble in
the triple bubble.

\subsection{Evolutions with different surface energies}

\subsubsection{3 phases}

As an example for non-equal surface energy densities for the various curves, we
repeat the simulation in Figure~\ref{fig:2d_db_plus_one}, but now weigh
curves 1 and 3 in the double bubble with $\sigma_1 = \sigma_3 = 2$, while
keeping the other two densities at unity. This now means that in contrast to
Figure~\ref{fig:2d_db_plus_one}, it makes energetically more sense to increase
the size of the single bubble, while shrinking the bubble that is surrounded by
the more expensive interfaces. See Figure~\ref{fig:2d_db_plus_one} for the
observed evolution.
\begin{figure}
\center
\includegraphics[angle=-90,width=0.18\textwidth]{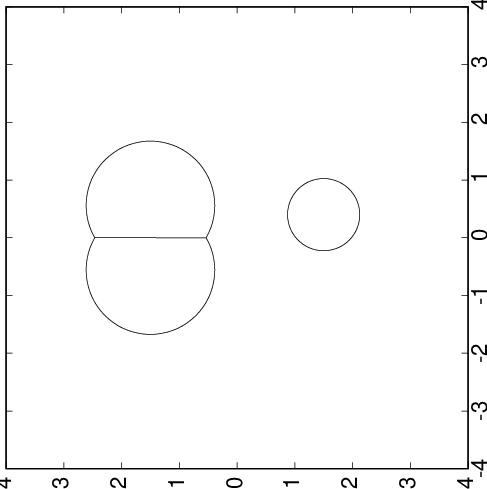}
\includegraphics[angle=-90,width=0.18\textwidth]{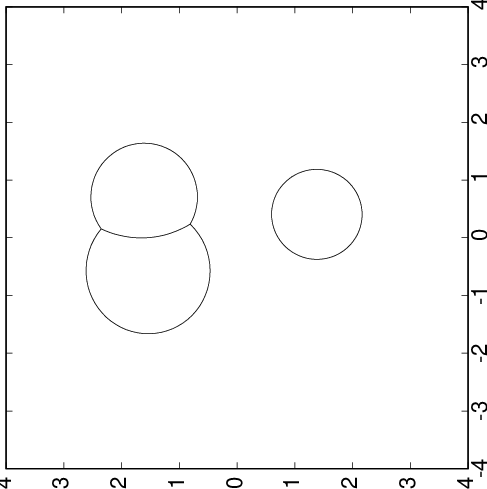}
\includegraphics[angle=-90,width=0.18\textwidth]{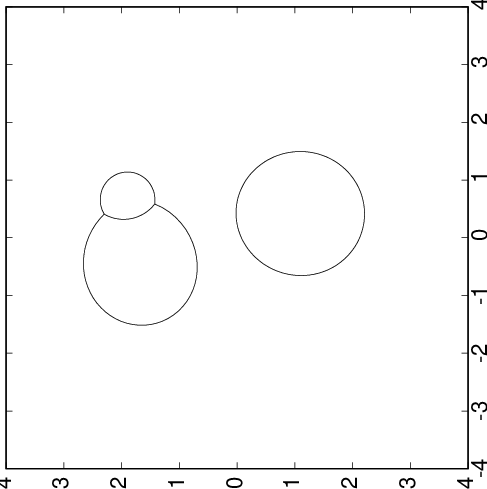}
\includegraphics[angle=-90,width=0.18\textwidth]{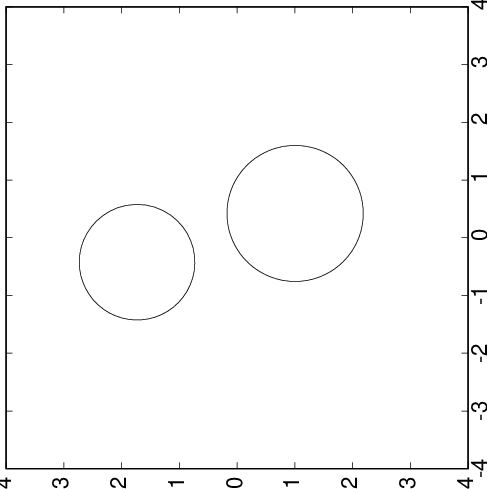}
\includegraphics[angle=-90,width=0.18\textwidth]{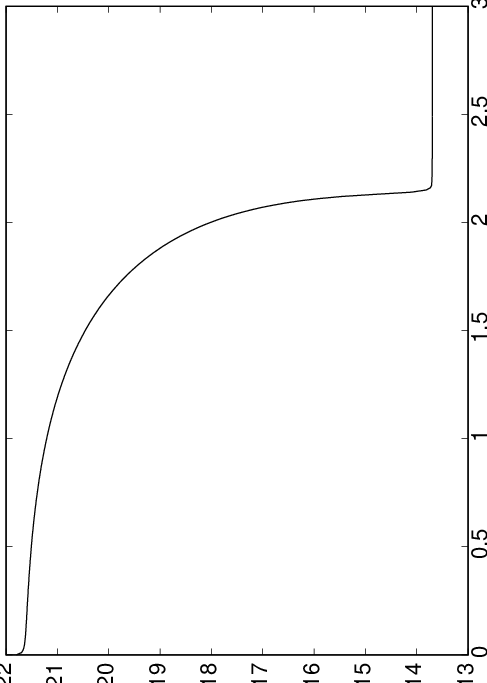}
\caption{The solution at times $t=0, 1, 2, 3$, and a plot of the discrete 
energy over time.
}
\label{fig:2d_iso2121_db_plus_one}
\end{figure}%

\subsubsection{4 phases}

\begin{figure}[h]
\center
\includegraphics[angle=-90,width=0.18\textwidth]{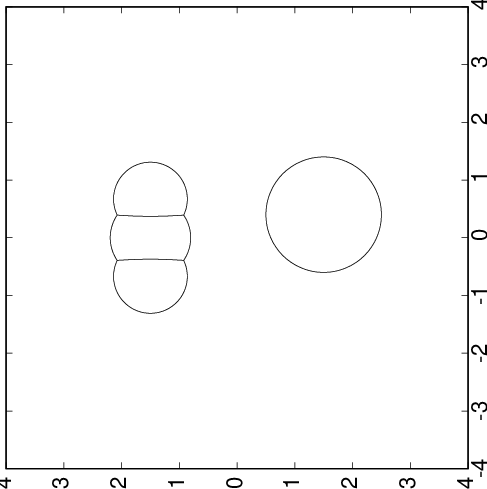}
\includegraphics[angle=-90,width=0.18\textwidth]{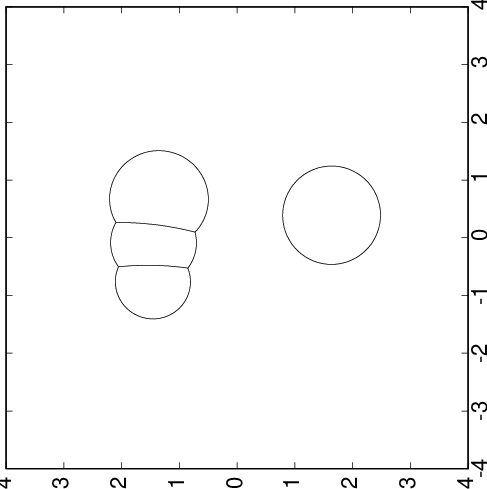}
\includegraphics[angle=-90,width=0.18\textwidth]{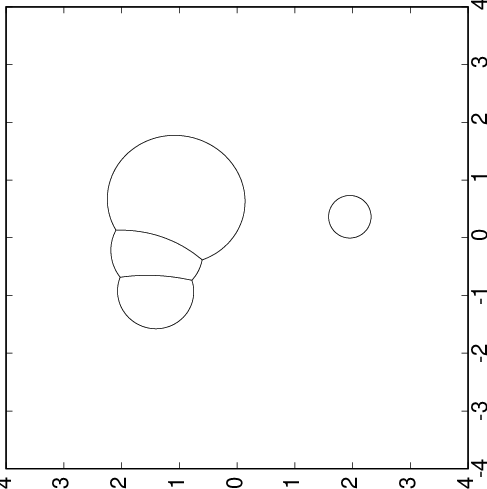}
\includegraphics[angle=-90,width=0.18\textwidth]{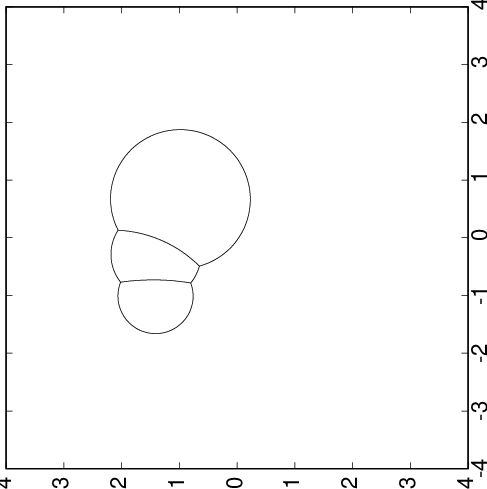}
\includegraphics[angle=-90,width=0.18\textwidth]{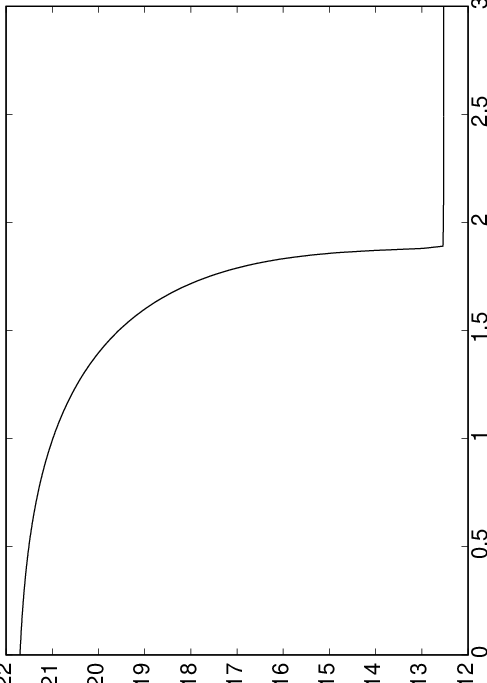}
\caption{The solution at times $t=0, 1, 1.8, 3$, and a plot of the discrete 
energy over time.
}
\label{fig:2d_iso1111112_tb_unstable_plus_bigone}
\end{figure}

Similarly, if we make the interface of the single circular bubble in
the initial data in Figure~\ref{fig:2d_tb_unstable_plus_bigone} more expensive,
it will no longer grow but shrink to a point. Setting the weight for the curve
to $\sigma_7 = 2$ leads to the evolution seen in 
Figure~\ref{fig:2d_iso1111112_tb_unstable_plus_bigone}.

\section{Acknowledgements}

A portion of this study was conducted during the research stay of the first author in Regensburg, Germany.
He is grateful to the administration of the University of Regensburg for the invitation and for the financial support
by the DFG Research Training Group 2339 \textit{IntComSin} - Project-ID 32182185.
He also appreciates the kindness and treatment received from all members in Research Training Group 2339.

\begin{appendix}
\renewcommand{\theequation}{\Alph{section}.\arabic{equation}}

\section{An exact solution for the three-phase Mullins--Sekerka flow}\label{subsec:DerivationExactSolution}
\label{sec:B}

In this appendix, we shall prove that \eqref{eq:ODE} 
is a solution to \eqref{eq:StrongForm_Multi} with $\sigma=1$.
Firstly, it directly follows from the definitions of $F$ and $\alpha$ in
\eqref{eq:ODE} that $F(R_2(t)) = \alpha(t) / R_2(t)$. Hence,
\eqref{eq:3circleRadiusFormula} immediately implies
that the normal velocity $V$ of $\curve{}{}(t)$ satisfies
    \begin{equation}\label{eq:3circleVelocityFormula}
        V = \alpha(t)
        \begin{cases}
            \frac1{R_1(t)} & \mbox{on}\ \ \curve{}{1}(t),\\
            -\frac1{R_2(t)} & \mbox{on}\ \ \curve{}{2}(t),\\
            \frac1{R_3(t)} & \mbox{on}\ \ \curve{}{3}(t).
        \end{cases}
    \end{equation}

{
    The first and fourth equations in $\eqref{eq:StrongForm_Multi}$ 
hold trivially since $\bv w$ as defined in \eqref{eq:3circleChemicalFormula} 
is constant in the two connected components of $\Omega_1(t)$, and harmonic
in $\Omega_2(t)$ and $\Omega_3(t)$.
    Let us confirm the Gibbs--Thomson law. 
We see from \eqref{eq:3circleChemicalFormula} that
$w_1 - w_2 = 1/R_1(t)$ on $\curve{}{1}(t)$, 
$w_3 - w_2 = -1/R_2(t)$ on $\curve{}{2}(t)$ and
$w_3 - w_1 = 1/R_3(t)$ on $\curve{}{3}(t)$.
    Thus, $\bv w$ satisfies the second condition in $\eqref{eq:StrongForm_Multi}$.
    We move on the motion law. A direct calculation shows 
    \begin{equation}\label{eq:3circleVelocityValidation}
        \begin{cases}
            -\jump{\nabla w_1}\cdot\normalnoindex{}{1} = \jump{\nabla w_2}\cdot\normalnoindex{}{1} = \frac{\alpha(t)}{R_1(t)}\ \ \mbox{on}\ \ \curve{}{1}(t),\\
            \jump{\nabla w_2}\cdot\normalnoindex{}{2} = -\jump{\nabla w_3}\cdot\normalnoindex{}{2} = -\frac{\alpha(t)}{R_2(t)}\ \ \mbox{on}\ \ \curve{}{2}(t),\\
            -\jump{\nabla w_1}\cdot\normalnoindex{}{3} = \jump{\nabla w_3}\cdot\normalnoindex{}{3} = \frac{\alpha(t)}{R_3(t)}\ \ \mbox{on}\ \ \curve{}{3}(t).
        \end{cases}
    \end{equation}
    Hence, the third condition of $\eqref{eq:StrongForm_Multi}$ is valid by $\eqref{eq:3circleVelocityFormula}$ and $\eqref{eq:3circleVelocityValidation}$. 
    Therefore, $\bv w = (w_1,w_2,w_3)^T$ given by $\eqref{eq:ODE}$ is an 
    exact solution of $\eqref{eq:StrongForm_Multi}$ with $\sigma=1$.
}
\end{appendix}

\end{document}